\documentclass[twoside]{amsart}

\usepackage[utf8]{inputenc}
\usepackage[T1]{fontenc}
\usepackage{t1enc}

\usepackage{amsmath,amssymb,amsfonts,amsthm,amscd,mathrsfs,mathtools}
\numberwithin{equation}{section}

\usepackage[mathscr]{eucal}

\usepackage{indentfirst}
\usepackage{relsize}
\usepackage[tablename=Table]{caption}
\usepackage[active]{srcltx}

\usepackage{graphicx}
\usepackage{graphics}
\usepackage{pict2e}
\usepackage{epic}
\usepackage{epstopdf}
\usepackage{tabularx}
\usepackage{booktabs,makecell,siunitx,eqparbox}

\usepackage[x11names,dvipsnames]{xcolor}

\usepackage[all]{xy}
\usepackage{tikz}
\usepackage{tikz-cd}
\usetikzlibrary{
  matrix,
  calc,
  arrows,
  arrows.meta,
  positioning,
  shapes,
  decorations.pathmorphing
}

\usepackage{ifpdf}

\usepackage{enumitem}
\setlist[enumerate]{label=\arabic*)}

\newcommand{\ba}{\begin{array}}
\newcommand{\ea}{\end{array}}

\newcommand{\op}[1]{{\mathcal O}_{\mathbb{P}^{#1}}}
\newcommand{\p}[1]{{\mathbb{P}^{#1}}}
\newcommand{\opn}{{\mathcal O}_{\mathbb{P}^{n}}}
\newcommand{\pn}{{\mathbb{P}^{n}}}
\newcommand{\tp}[1]{{\rm T}\mathbb{P}^{#1}}
\newcommand{\tpn}{T{\mathbb{P}^{n}}}

\newcommand{\sI}{\mathscr{I}}

\newcommand{\dV}{{\mathrm{dV}}}

\DeclareMathOperator{\sing}{Sing}
\DeclareMathOperator{\Supp}{Supp}
\DeclareMathOperator{\im}{im}
\DeclareMathOperator{\rk}{rk}

\DeclareMathOperator{\Hom}{Hom}
\DeclareMathOperator{\ch}{ch}

\DeclareMathOperator{\codim}{codim}
\DeclareMathOperator{\Ext}{Ext}

\DeclareMathOperator{\id}{id}

\DeclareMathOperator{\coker}{coker}
\DeclareMathOperator{\rad}{rad}

\DeclareMathOperator{\depth}{depth}

\newcommand{\sln}{\mathfrak{sl}}
\newcommand{\gl}{\mathfrak{gl}}
\newcommand{\mf}[1]{\mathfrak{#1}}
\DeclareMathOperator{\gst}{gst}
\DeclarePairedDelimiter{\bkt}{[}{]}
\DeclareMathOperator{\SL}{SL}

\def\PP{{\mathbb{P}}}
\def\OO{\mathcal{O}}
\def\CC{{\mathbb{C}}}

\newcommand{\QQ}{\mathbb{Q}}

\newcommand{\sF}{\mathscr{F}}
\newcommand{\sG}{\mathscr{G}}
\newcommand{\sH}{\mathscr{H}}

\newcommand{\tF}{{\rm T}\sF}
\newcommand{\tG}{{\rm T}\sG}

\newcommand{\nF}{{\rm N}\sF}

\newcommand{\nH}{{\rm N}\sH}

\newcommand{\Fol}{\mathbb{F}{\rm ol}}

\newcommand{\tX}{{\rm T}X}

\newcommand{\lra}{\longrightarrow }
\newcommand{\aut}{\mathfrak{aut}}
\newcommand{\fix}{\mathfrak{fix}}
\newcommand{\afix}{\mathfrak{afix}}
\DeclareMathOperator{\Aut}{Aut}
\DeclareMathOperator{\PGL}{PGL}
\newcommand{\lieder}{\mathcal{L}}
\newcommand{\extp}{{\textstyle \bigwedge}} 
 
\newcommand{\del}[1]{\frac{\partial}{\partial #1}}

\newcommand{\lpb}[3]{{\rm LPB }(#1,#2,#3)}

\newcommand{\logg}[3]{{\rm Log}(#1;#2,#3)}
\newcommand{\exc}[3]{{\rm ALie}(#1,#2;#3)}

\newcommand{\cont}{\mathbin{\lrcorner}}
\newcommand{\fp}{\mathfrak{p}}
\newcommand{\NC}{\mathcal{NC}}

\tikzset{
  >=latex,
  vertex/.style={circle,draw,inner sep=2pt},
  active/.style={thick},
  inactive/.style={dashed}
}

\tikzset{
  vertex/.style={circle,draw,inner sep=1pt,fill=white},
  edgeActive/.style={line width=0.9pt,color=RoyalBlue3},
  edgeInactive/.style={line width=0.4pt,color=Gray},
}

\newtheorem{theorem}{Theorem}[section]
\newtheorem{proposition}[theorem]{Proposition}
\newtheorem{corollary}[theorem]{Corollary}
\newtheorem{lemma}[theorem]{Lemma}

\newtheorem{innerthmintro}{Theorem}
\newenvironment{thmintro}[1]
{\renewcommand\theinnerthmintro{#1}\innerthmintro}
{\endinnerthmintro}

\newtheorem{innerpropintro}{Proposition}
\newenvironment{propintro}[1]
{\renewcommand\theinnerpropintro{#1}\innerpropintro}
{\endinnerpropintro}

\theoremstyle{definition}

\newtheorem{example}{Example}[section]
\newtheorem{remark}{Remark}[section]

\usepackage{hyperref}
\hypersetup{
  colorlinks=true,
  linkcolor=blue,   
  filecolor=magenta,
  urlcolor=cyan
}

\date{}

\begin{document}

\title{ Holomorphic foliations of degree two and arbitrary dimension }

\author[M. Corr\^ea]{ Maur\'icio Corr\^ea }
\address[M. Corr\^ea]{Dipartimento di Matematica, Universit\`a degli Studi di Bari,  Via E. Orabona 4, I-70125, Bari, Italy
}
\email{mauricio.correa.mat@gmail.com, mauricio.barros@uniba.it} 

\author[A. Muniz]{Alan Muniz }
\address[A. Muniz]{Universidade Estadual de Campinas (UNICAMP) \\ Instituto de Matemática, Estatística e Computação Científica (IMECC) \\ Departamento de Matem\'atica \\
Rua S\'ergio Buarque de Holanda, 651\\ 13083-859 Campinas-SP, Brasil}
\address{Departamento de Matem\'atica \\ Centro de Ci\^encias Exatas e da Natureza \\ Universidade Federal de Pernambuco \\ Recife - PE, CEP 50740-560, Brasil }
\email{alannmuniz@gmail.com, alan.nmuniz@ufpe.br }



\begin{abstract}
We give a complete classification of degree-two holomorphic foliations on projective spaces. Up to linear pullback, these are logarithmic foliations, linear pullback foliations, affine-Lie foliations induced by an action of the affine Lie algebra, or foliations of new types whose tangent sheaves are linear pullbacks of either $\tp3(-1)$ or a null-correlation bundle. We call the foliations arising from null-correlation bundles null-correlation foliations. For two-dimensional foliations, we obtain a complete decomposition of the space of foliations into irreducible components. The null-correlation component is generically reduced for $n\geq8$ and generically nonreduced for $n=7$. 
\end{abstract}

\maketitle


\section{Introduction}

One of the central problems in the global theory of holomorphic foliations is the description and classification of the irreducible components of the space of foliations on projective space, at least in low degrees. In degree zero, every foliation is induced by a linear rational map, and the corresponding space is a Grassmannian \cite[Th\'eor\`eme~3.8]{CD2}. In degree one, Loray, Pereira and Touzet \cite{LPT2}, extending results of Jouanolou \cite{Jou}, proved that there are two irreducible components, one formed by linear pullbacks of foliations by curves and the other by rational foliations. For degree two, Cerveau and Lins Neto \cite{CLN} classified codimension-one foliations and obtained six irreducible components. Since then, several works have investigated the classification in dimension $n-1$ and higher degrees \cite{CLNdeg3, AVdeg4, LPT3}, most notably in \cite{CLPdeg3}, where the authors present a partial classification for degree $3$. In a different direction, the present work focuses on degree-$2$ foliations of higher codimension. Our first main result is a complete description of the space of foliations of dimension and degree equal to two.

\begin{thmintro}{A}\label{thm:mainA}
The space $\Fol(2,2,n)$ of two-dimensional foliations of degree two on $\pn$, $n\geq 3$, decomposes into the following irreducible components. 
\begin{enumerate}
    \item Four logarithmic components $\logg{1^{n-3},2,2}{2}{n}$, $\logg{1^{n-2},3}{2}{n}$, $\logg{1^{n-1},2}{2}{n}$, and $\logg{1^{n+1}}{2}{n}$;
    \item The linear pullback component $\lpb{2}{2}{n}$;
    \item Affine-Lie components $\exc{r_1,\dots,r_m}{2}{n}$, where $r_1+\dots+r_m=n+1$ is a partition different from the ones in the logarithmic case;
    \item The null-correlation component $\NC_n$, if $n\geq7$.
\end{enumerate}
\end{thmintro}

The notation is explained in \S\ref{sectionfol}. Linear pullback components go back to \cite{CLNpullback}; for the version needed here see \cite{CP1}. Logarithmic components are described in \cite{Omegar, CLNlog,CPV1}. We call Affine-Lie foliations those defined by actions of the Lie algebra $\mf{aff(\CC)}$ of affine transformations of the line. They generalize the exceptional component from \cite{CLN}. We refer to Pereira and dos Santos \cite{PdS} for a comprehensive description of such foliations. The null-correlation component $\NC_n$ is new. The tangent sheaf of a null-correlation foliation is the reflexive pullback of a null-correlation bundle under a linear projection to $\p3$. The closure of a general leaf is a projective plane. 
 
The component $\NC_n$ is the closure of a projective orbit (see \S\ref{sec:nullcorr-component}). We use curl-free bivectors to study first-order deformations along this orbit. For $n\geq8$, every first-order deformation of a general null-correlation foliation is induced by a projective change of coordinates. For $n=7$, the orbit has dimension $40$ and the Zariski tangent space to $\Fol(2,2,7)$ at a point of the orbit has dimension $44$.

We prove Theorem \ref{thm:mainA} in \S\ref{sec:class-dim2}. The proof is organized according to the stability of the tangent sheaf. If $\tF$ is stable, we study the reflexive hull of the cokernel of the lift $\tF\hookrightarrow\opn(1)^{\oplus n+1}$. Its splitting gives rational foliations (Lemma \ref{lem:split-to-rational}), while the nonsplit case leads to the null-correlation foliations. If $\tF$ is strictly semistable, we study the subfoliations generated by tangent global vector fields and use the associated Lie algebras and algebraic group actions to obtain logarithmic and affine-Lie cases. If $\tF$ is unstable, a nonzero section of $\tF(-1)$ gives a linear projection under which $\sF$ is the pullback of a degree-two foliation by curves.

For the higher-dimensional case we have the following result.

\begin{thmintro}{B}\label{thm:mainC-intro}
Let $\sF$ be a foliation of degree two and dimension $k>2$ on $\pn$, and put $q=n-k$. Then one of the following holds.
\begin{enumerate}
    \item $\sF$ belongs to a logarithmic component;
    \item $\sF$ belongs to an Affine-Lie component;
    \item $\sF$ belongs to the Linear Pullback component;
    \item $\sF$ is the linear pullback of a null-correlation foliation on $\p{n-k+2}$, with $n\geq k+5$;
    \item $\sF$ is the linear pullback of a three-dimensional degree-two foliation $\sG$ on $\p{n-k+3}$ whose tangent sheaf is the reflexive pullback of $\tp3(-1)$ under a linear projection $\rho\colon\p{n-k+3}\dashrightarrow\p3$.
\end{enumerate}
\end{thmintro}

To prove Theorem \ref{thm:mainC-intro}, we note that $\sF$ is the linear pullback (possibly by the identity) of a foliation $\sG$ such that $H^0(\tG(-1)) = 0$. If $\dim \sG = 1$, then $\sF$ is a Linear Pullback foliation. If $\dim \sG = 2$, then the result follows from Theorem \ref{thm:mainA}. If $\dim \sG >2$, then we have either (1) or (5); see Proposition \ref{prop:ra=r-1=log}.

Finally, we remark that the number of irreducible components of $\Fol(2,k,n)$ grows with the number of partitions of $n-k+3$, as already noticed in \cite[Corollary D]{PdS}. 

In proving Theorem~\ref{thm:mainA}, we obtained the following result, which we prove in Section~\ref{sec:subfol}. We observe that this result was also obtained by da Costa \cite{C} with independent techniques.

\begin{propintro}{D}\label{prop:B}
Let $\sF$ be a two-dimensional foliation of degree $d$ on $\pn$ such that $h^0(\tF)>0$. If $\tF$ is locally free, then there is $a\in\{0,1\}$ such that 
\[
\tF\cong\opn(a)\oplus\opn(2-d-a).
\]
\end{propintro}

The paper is organized as follows. Section~\ref{sectionfol} fixes the notation and recalls the standard degree-two families. Section~\ref{sec:subfol} studies two-dimensional foliations tangent to global vector fields. Section~\ref{sec:stable-dim2} treats stable tangent sheaves. Section~\ref{sec:nullcorr-component} studies the null-correlation component in dimensions $n\geq8$ and $n=7$. Section~\ref{sec:class-dim2} completes the classification in dimension two and gives the reduction in higher dimension.

\subsection*{Usage of Large Language Models}
During the development and preparation of this article, the authors used ChatGPT (OpenAI; GPT-6.1 Sol) as an interactive tool to explore proof strategies, search the literature, and assist in revision. The main contribution of the LLM was to point out the reference \cite{BeRo} that helped control the cokernel of $\tF \hookrightarrow \opn(1)^{\oplus n+1}$. It also helped with minor corrections and computations.

\subsection*{Acknowledgments}
We thank Jorge Vit\'orio Pereira for his interest in this work and for numerous insightful conversations. MC is partially supported by the Universit\`a degli Studi di Bari and by the PRIN 2022MWPMAB- ``Interactions between Geometric Structures and Function Theories'' and is a member of INdAM-GNSAGA. AM was partially supported by INCTmat/MCT/Brazil, CNPq grant number 160934/2022-2, and also by a CNPq research productivity fellowship, grant number 305599/2026-7.

\section{Holomorphic foliations}\label{sectionfol}
In this section, we review the basic concepts needed throughout the work. Let $X$ be a smooth complex projective variety. A \emph{singular holomorphic distribution} $\sF$ on $X$ is the data of a saturated subsheaf $\tF$ of the tangent sheaf $\tX$ of $X$. These sheaves give a short exact sequence
\begin{equation}\label{eq:dist-def}
    0  \lra  \tF  \lra  \tX  \lra  \nF  \lra  0,
\end{equation}
where $\nF$ is a torsion-free sheaf, called the \emph{normal sheaf} of $\sF$, and $\tF$ is a reflexive sheaf, called the \emph{tangent sheaf} of $\sF$. The dimension of $\sF$ is the generic rank of $\tF$. If $\tF$ is involutive, i.e. $[\tF, \tF] \subset \tF$, then $\sF$ is a foliation; this means that through a general point in $X$ passes a unique connected immersed (analytic) subvariety $L \subset X$ such that ${\rm T}L = \tF|_L$, called a \emph{leaf}. A foliation $\sF$ is \emph{algebraically integrable} if its leaves are embedded algebraic varieties. 

Taking the top exterior power of the quotient map in \eqref{eq:dist-def}, we obtain a morphism $\omega_{\sF}\colon \extp^p\tX \to \det(\nF)$, where $p = \rk \nF$, whose image is $\det(\nF)$ tensored by the ideal sheaf $\sI_Z$ of the vanishing locus of $\omega_{\sF}$ seen as a global section of $\Omega_X^p\otimes \det(\nF)$. We call $Z$ the \emph{singular scheme} of $\sF$, and denote it by $\sing(\sF)$. Since $\tF$ is saturated, $\codim(\sing(\sF))\geq 2$.

Outside the singular scheme, the twisted $p$-form $\omega_\sF\in H^0(\Omega_X^p\otimes \det(\nF))$ is locally decomposable: for any $p \in X \setminus \sing(\sF)$ there exists an open neighborhood $U\ni p$ and $1$-forms $\eta_1, \dots, \eta_p \in \Omega_X^1(U)$ such that 
$
\omega_\sF|_U = \eta_1\wedge \dots \wedge \eta_p.
$
Moreover, $\sF$ is a foliation if and only if
$
d\eta_j \wedge \omega_\sF|_U = 0, \, \text{for} \, j= 1, \dots, p.
$

\subsection{Foliations on projective spaces}\label{sect:folsonPn}
A $k$-dimensional distribution $\sF$ on $\pn$ corresponds to a twisted $(n-k)$-form $\omega \in H^0(\Omega_{\pn}^{n-k}(d+n-k+1))$, where $d\geq 0$ is called the \emph{degree} of $\sF$. The Euler sequence allows us to interpret $\omega$ as a homogeneous polynomial $(n-k)$-form
\[
\omega=\sum_{0\leq i_1<\cdots<i_{n-k}\leq n}
A_{i_1\dots i_{n-k}}\,
 dx_{i_1}\wedge\cdots\wedge dx_{i_{n-k}},
\]
where the coefficients $A_{i_1\dots i_{n-k}}$ are homogeneous of degree $d+1$, and the contraction with the radial vector field $\rad = x_0\del{ x_0} + \cdots + x_n\del{ x_n}$ vanishes: ${\rad} \cont \omega = 0$. The local decomposability translates to the Pl\"{u}cker conditions  
\begin{equation}\label{eq:lds}
    (v \cont \omega) \wedge \omega = 0 \quad \text{ for every }v
 \in \wedge^{n-k-1}(\mathbb C^{n+1}) ,
\end{equation}
and the integrability condition becomes
\begin{equation}\label{eq:integrability}
     (v \cont \omega) \wedge d\omega = 0 \quad \text{ for every }v
 \in \wedge^{n-k-1}(\mathbb C^{n+1}) .
\end{equation}
We refer to \cite{Jou, CPV1} for details.

\subsubsection{Integrability of bivectors} In the case of two-dimensional foliations, $H^0(\Omega_\pn^{n-2}(d+n-1)) \cong H^0(\extp^2\tpn(d-2))$. The latter can be seen as homogeneous bivectors modulo multiples of the radial vector field. Given $\pi \in H^0(\extp^2\tpn(d-2))$ the conditions \eqref{eq:lds} and \eqref{eq:integrability} become
\begin{equation}\label{eq:int-bivector}
\pi \wedge \pi \wedge \rad  = 0 \quad \text{and} \quad \bkt{\pi,\pi}\wedge \rad  = (d-2) \pi \wedge \pi,
\end{equation}
where $\bkt{\, ,\,}$ is the Schouten bracket. Fix the volume form $\dV = dx_0\wedge \dots \wedge dx_n$. For a $p$-vector $\xi$, let $D(\xi)$ be the $(p-1)$-vector such that $d( \xi \cont \dV) = D(\xi)\cont \dV$. This is called the \emph{curl} of $\xi$. If $p=1$, this recovers the divergence of a vector field. By Koszul's Theorem \cite[Theorem 2.6.7]{DN-poisson}, if $\xi$ is a $p$-vector and $\zeta$ is a $q$-vector, then 
\begin{equation}\label{eq:curl-bracket}
D(\xi\wedge \zeta) = (-1)^qD(\xi)\wedge\zeta  + \xi \wedge D(\zeta) + (-1)^{q}\bkt{\xi,\zeta}. 
\end{equation}
Given a bivector $\pi$ of total degree $d-2$ as above, consider $\pi_0 = \pi - \frac{1}{(d+n-1)}D(\pi) \wedge \rad$. Then $\pi_0$ and $\pi$ define the same element in $H^0(\extp^2\tpn(d-2))$; the formula above implies $D(\pi_0) = 0$. Hence, we can always choose a curl-free bivector to define the foliation. 
In this case, \eqref{eq:int-bivector} becomes the following.

\begin{lemma}\label{lem:curl-free-integ}
Let $\pi$ be a curl-free homogeneous bivector of total degree $d-2$. Then $\pi$ is integrable if and only if $\pi\wedge \pi = 0$. 
\end{lemma}

\begin{proof}
Since $D(\pi) = 0$, we have $D(\pi\wedge\pi) = \bkt{\pi,\pi}$. In particular, $D(\bkt{\pi,\pi}) = 0$. Suppose $\pi$ is integrable. By \eqref{eq:int-bivector} 
\[
0 = D(\pi \wedge \pi \wedge \rad)  = -D(\pi\wedge\pi)\wedge \rad +(n+1+2d-4) \pi \wedge \pi
\]
Hence $\bkt{\pi,\pi}\wedge \rad = (n+2d-3) \pi \wedge \pi$. By \eqref{eq:int-bivector} again, $(n+d-1) \pi \wedge \pi = 0$. Since $n+d-1\neq 0$, $\pi \wedge \pi = 0$. Conversely, if $\pi \wedge \pi = 0$, then $\bkt{\pi,\pi} = 0$ and $\pi$ is integrable.
\end{proof}

\subsubsection{Lift of vector fields and global decomposability}\label{sec:lift-decomp}

Let $S = \CC[x_0,\dots, x_n]$ be the polynomial ring, $T_S = Der_\CC(S) \cong S(1)^{n+1}$ the module of derivations and $\Omega_S^1 = {T_S}^\vee$ the module of K\"ahler differentials. Given a homogeneous polynomial $p$-form $\eta$ on $\CC^{n+1}$ of total degree $d+p+1$, one defines the contraction map
\[
\cont \eta \colon T_S \lra \Omega_S^{p-1}(d+p+1),
\]
whose kernel is the module of vector fields $v \in T_S$ that annihilate $\eta$. Note that $\eta$ is LDS if and only if $\ker \cont \eta$ is a submodule of generic rank $n+1-p$.

\begin{lemma}\label{lem:decomp-free}
Let $\eta \in \Omega_S^p$ be an LDS $p$-form without zeros in codimension one and let $M = \im \cont \eta$. Then $\eta$ is globally decomposable if and only if $M^{\vee\vee}$ is a free $S$-module.
\end{lemma}

\begin{proof}
If $\eta =  \theta_1 \wedge \dots \wedge \theta_p$ then $M^{\vee\vee}$ is the free $S$-module generated by the $(p-1)$-fold products  $\eta_j = \theta_1 \wedge \dots \wedge \widehat{\theta_j}\wedge \dots \wedge\theta_p$. Indeed, if $N = \oplus_jS\eta_j$, then $M\subset N$ and $M_\fp = N_\fp$ for any prime $\fp$ of codimension one, since $\codim \sing(\eta) \geq 2$. Then  $M^{\vee\vee} = \bigcap_{\codim \fp = 1}M_\fp =  N$, see \cite[\href{https://stacks.math.columbia.edu/tag/0AUY}{Tag 0AUY}]{stacks-project}. 

Conversely, suppose $M^{\vee\vee} \cong \bigoplus_j S(a_j)$ is free. Then the image of the dual morphism
\[
(\cont \eta)^\vee \colon M^\vee = \bigoplus_j S(-a_j)  \lra \Omega_S^1
\]
defines independent $1$-forms $\theta_1, \dots, \theta_p$ such that $\theta_i \wedge \eta = 0$. Then Saito's Division Lemma implies that $f\eta = \theta_1 \wedge \dots \wedge \theta_p$ for some polynomial $f$. On the other hand, $\cont \eta$ factors through $\psi \colon T_S \lra M^{\vee\vee}$ given by $(\cont \theta_1, \dots, \cont \theta_p)$. For any prime $\fp$ of codimension one, $\psi_\fp$ is surjective since $M_\fp = M_\fp ^{\vee\vee}$. Therefore, $\psi$ does not drop rank in codimension one, i.e.,  $\theta_1 \wedge \dots \wedge \theta_p$ does not vanish in codimension one. Since $\codim \sing(\eta) \geq 2$, $f$ is a constant.
\end{proof}

We are interested in the case $\eta = d\omega$ where $\omega$ defines a foliation on $\pn$. Indeed, since $\omega$ is integrable, $d\omega$ is LDS \cite[Remark 1.2.1.]{Medeiros}. Hence Lemma \ref{lem:decomp-free} hints that if the corresponding $M^{\vee\vee}$ is free, then $\sF$ is a rational foliation. This is true in degree two if we suppose $H^0(\tF^\vee) = 0$, see Lemma \ref{lem:split-to-rational} below. The following is a trivial fact that lacks a precise reference.

\begin{lemma}\label{lem:lift-map}
Let $\sF$ be a foliation on $\pn$. Then $\tF$ injects into $\opn(1)^{\oplus n+1}$ as the sheaf associated with $\ker \cont d\omega$. 
\end{lemma}

\begin{proof}
Define $F = \widetilde{\ker \cont d\omega}$. Then $F$ is reflexive of rank $n-p$ with an injection $F\hookrightarrow  \opn(1)^{\oplus n+1}$. Then the inclusion $\ker d\omega \subset \ker \omega$ yields $\tF \cong (F + \opn\cdot\rad)/(\opn\cdot\rad) \cong F$.
\end{proof}

\subsection{Irreducible components of the space of foliations}
Given integers $d$, $k$, and $n$, the conditions \eqref{eq:lds} and \eqref{eq:integrability} define the following subvariety of $\PP H^0(\Omega_{\pn}^{n-k}(n-k+d+1))$:
\begin{equation}\label{eq:def-Fol}
    \Fol(d,k,n) =  \Big\{[\omega]   \mid \omega \text{ is integrable, and } \codim( \sing \omega) \geq 2 \Big\},
\end{equation}
the space of $k$-dimensional foliations of degree $d$ on $\pn$. For $d \leq 1$ the irreducible components of $\Fol(d,k,n)$ are well-known. A degree-$0$ foliation is defined by a linear rational map $\pn \dashrightarrow \p{n-k}$, see \cite[ Th\'eor\`eme 3.8]{CD2}. In particular,
\[
\Fol(0,k,n) \cong G(k,n+1),
\]
the Grassmannian of $k$-dimensional vector subspaces of $\CC^{n+1}$. For degree one, there is the following classification. 

\begin{theorem}\cite[ Theorem 6.2]{LPT2}\label{thm:fols-deg1}
If $\sF$ is a foliation of degree $1$ and codimension $q$ on $\pn$, then we are in one of the following cases:
\begin{enumerate}
\item $\sF$ is defined by a dominant rational map $\rho \colon \pn \dashrightarrow \PP(1^{q},2)$, with an irreducible general fiber determined by $q$ linear polynomials and one quadratic form; or
\item $\sF$ is the linear pullback of a foliation induced by a global holomorphic vector field on $\mathbb{P}^{q+1}$. 
\end{enumerate}
\end{theorem}
Hence, $\Fol(1,k,n)$ has two irreducible components. In degree $2$, the codimension-one case was classified by Cerveau and Lins Neto \cite{CLN}. 

\begin{theorem}\cite[Theorem A]{CLN}
The space of codimension-one foliations of degree $2$ in $\pn$, $n\geq 3$, has $6$ irreducible components: 
\begin{itemize}
    \item ${S(2,n)}$ composed of linear pullbacks of degree-$2$ foliations on $\p2$;
    \item ${R(2,2)}$ given by pencils of quadrics;
    \item ${R(1,3)}$ given by pencils of cubics with a triple hyperplane;
    \item ${L(1,1,1,1)}$ whose general member is defined by a logarithmic $1$-form with poles in $4$ hyperplanes;
    \item ${L(1,1,2)}$ whose general member is defined by a logarithmic $1$-form with poles in two hyperplanes and a quadric;
    \item ${E(n)}$ the \emph{exceptional component} whose general member is defined by the linear pullback of an algebraic action of the Lie algebra of affine transformations $\mathfrak{aff}(\CC)$.
\end{itemize}
\end{theorem}

Below, we detail each one of these components. We describe their generalizations to any codimension. 

\begin{example}[Linear pullback foliations]\label{ex:lpb-fols}
Given a linear projection $\rho \colon \pn \dashrightarrow \p{n-k+1}$ and a degree-$d$ foliation by curves $\sG$ on $\p{n-k+1}$, the pullback $\sF=\rho^*\sG$ is a degree-$d$ foliation of dimension $k$. Varying $\rho$ and $\sG$, we get an irreducible component $\lpb{d}{k}{n} \subset \Fol(d,k,n)$. For a proof, see \cite[Corollary 5.1]{CP1}.
\end{example}

\begin{example}[Logarithmic foliations]\label{ex:log-fols}
A logarithmic foliation on $\pn$ of dimension $k$ is defined by a closed rational $(n-k)$-form
\[
\begin{split}
\omega={}&\sum_{0\leq i_1<\cdots<i_{n-k}\leq r}\lambda_{i_1,\dots,i_{n-k}}
\frac{dF_{i_1}}{F_{i_1}}\wedge\cdots\wedge
\frac{dF_{i_{n-k}}}{F_{i_{n-k}}},
\end{split}
\]
where $F_0,\dots,F_r$ are homogeneous polynomials of degrees $d_0,\dots,d_r$ and $r\geq n-k$. The coefficients $\lambda_{i_1,\dots,i_{n-k}}\in\CC$ satisfy the Pl\"ucker conditions, and $\omega$ is assumed to have no divisorial zeros. By \cite[Theorem 1.17]{CLNlog}, if $r\geq n-k+1$ these foliations form an irreducible component 
\[
\logg{d_0, \dots, d_r}{k}{n}\subset   \Fol( d, k, n), 
\]
where $d = -n+k-1 + \sum_{j=0}^{r} d_j$. When $r=n-k$, these are also called rational foliations; see \cite{CPV1}. 
\end{example}

\begin{remark}\label{rem:log(1,2)}
In degree $2$, the equation $2=-n+k-1+\sum_{j=0}^{r}d_j$, $r\geq n-k$, has the following solutions
\[
(1^{n-k+1},2), \quad  (1^{n-k+3}), \quad (1^{n-k-1}, 2,2), \quad \text{or}\quad  (1^{n-k}, 3).
\]
Hence $\Fol(2,k,n)$ has four logarithmic components. The last two are composed of rational foliations. A general foliation of type $(1^{n-k+1},2)$ is the pullback by a rational map $\rho\colon\pn\dashrightarrow\PP(1^{n-k+1},2)$
of a foliation induced by a global vector field on the weighted
projective space, and conversely. 
If $k>2$, a logarithmic foliation of type $(1^{n-k+3})$ is a linear
pullback of a two-dimensional foliation of the same type, since it is
defined, up to a linear change of coordinates, by fewer than $n+1$
variables. If $k=2$, a general foliation of type $(1^{n+1})$ is defined
by two commuting global vector fields on $\pn$, and hence by an action
of $\CC^2\hookrightarrow H^0(\tpn)$.
\end{remark}

\subsection{Foliations associated with a Lie algebra action}\label{sec:lie-alg-fols}
A $k$-dimensional Lie subalgebra $\mf{g} \subset \mf{sl}_{n+1} =  H^0(\tpn)$ defines an element $\omega(\mf{g}) \in H^0(\Omega_{\pn}^{n-k}(n+1))$. If $\omega(\mf{g})$ is nonzero and does not vanish in codimension one, it defines a foliation of degree equal to its dimension, $d=k$. We refer to \cite{CD2,CP1,PdS} for generalities and details about these foliations. Here we focus on $2$-dimensional Lie algebras. In \cite{PdS}, Pereira and dos Santos describe the scheme of Lie subalgebras and the irreducible components of the space of foliations they define.

Consider the algebraic variety of $2$-dimensional Lie subalgebras of $\sln_{n+1}$: 
\[
{\rm SLie }(2, \sln_{n+1}) = \{\,V \subset \sln_{n+1} \mid [V,V]\subset V, \, \dim V = 2 \,\}.
\]
It is a closed subvariety of the Grassmannian $G(2, \sln_{n+1})$ defined by the equations imposed by the Lie bracket. ${\rm SLie }(2, \sln_{n+1})$ has several irreducible components, one $\mathcal{C}$ corresponding to commutative Lie algebras and another for each nontrivial integer partition of $n+1$:
\[
{\rm SLie }(2, \sln_{n+1}) = \mathcal{C} \cup \bigcup_{r_1+\cdots +r_m = n+1}\mathcal{S}_{(r_1, \dots, r_m)}.
\]

Consider the map $\mathcal{F} \colon {\rm SLie }(2, \sln_{n+1})   \lra  \Fol(2,2,n)$ associating with a Lie subalgebra $\mathfrak{g}$ the foliation $\mathcal{F}(\mathfrak{g})$ it defines, whenever it makes sense. 
We saw in Remark \ref{rem:log(1,2)} that $\mathcal{F}(\mathcal{C}) = {\rm Log}(1^{n+1};2,n)$ is the logarithmic component. For the other components, there is the following result, which is a restatement of \cite[Theorem~3.6]{PdS}. They prove the statement for partitions having a part equal to $5$; the same argument applies to the remaining partitions, with the exceptional cases recorded below. 

\begin{theorem}\cite[Theorem~3.6]{PdS}\label{thm:lie-comps}
$\exc{r_1, \dots, r_m}{2}{n} := \mathcal{F}(\mathcal{S}_{(r_1, \dots, r_m)})$ is an irreducible component of $\Fol(2,2,n)$ unless $(r_1, \dots, r_m)$ is one of the following partitions:
\[
\{ (3,1^{n-2}), (2,2,1^{n-3}), (2, 1^{n-1}), (1^{n+1}) \}.
\]
Moreover, $\mathcal{F}$ is not defined on $\mathcal{S}_{(1^{n+1})}$ and $\mathcal{S}_{(2, 1^{n-1})}$.
\end{theorem}

\begin{remark}\label{Alg: log(1,2)} 
The partition $(4)$ gives the classical exceptional component on $\p3$ \cite{CLN}. The partition $(n+1)$ for $n\geq3$ appears in \cite[Example~6.8]{CP1}. These are all rigid, i.e., the component $\exc{n+1}{2}{n}$ is the orbit closure of a general foliation in it. We also remark that $\mathcal{F}(\mathcal{S}_{(3,1^{n-2})})$  and $\mathcal{F}(\mathcal{S}_{(2,2,1^{n-3})})$ do not fill irreducible components of $\Fol(2,2,n)$ and both fall into the logarithmic component $\logg{1^{n-1},2}{2}{n}$, see Proposition \ref{prop:class-om-nilp} below.
\end{remark}

\begin{remark}\label{rem:exceptional}
By \cite[Corollary~5.1]{CP1}, if $k\geq3$ and $r_1 + \dots + r_m = n-k+3$, there are irreducible components $\exc{r_1, \dots, r_m}{k}{n} \subset \Fol(2,k,n)$ whose general element is a linear pullback of a foliation in $\exc{r_1, \dots, r_m}{2}{n-k+2}$.  We call these the \emph{Affine-Lie foliations}.   
\end{remark}

\section{Study of subfoliations}\label{sec:subfol}
In this section, we study subfoliations $\sG \subset \sF$ of a fixed foliation $\sF$, with special attention to subfoliations defined by actions of Lie algebras. A Lie algebra naturally attached to a foliation $\sF$ on $\pn$ is given by its infinitesimal automorphisms
\[
\aut(\sF) = \{\, v \in H^0(\tpn) \mid \exists c\in \CC \text{ s.t. } \lieder_v\omega =c\,\omega \,\},
\]
where $\lieder_v$ stands for the Lie derivative in the direction of $v$ and $\omega$ is a twisted differential form defining $\sF$. $\aut(\sF)$ is the Lie algebra of the closed algebraic subgroup $\Aut(\sF) \subset \PGL(n+1,\CC)$ that preserves $\sF$. The global vector fields tangent to $\sF$ belong to $H^0(\tF)$ and define a Lie subalgebra of $\aut(\sF)$, usually denoted by $\fix(\sF)$. The exponential of $\fix(\sF)$ is not necessarily a closed subgroup. Let $\afix(\sF)$ denote the Lie algebra of the smallest algebraic subgroup whose Lie algebra contains $\fix(\sF)$. Note that
\[
\fix(\sF) \subset \afix(\sF) \subset \aut(\sF).
\]
The computation of $\afix(\sF)$ may be carried out as in \cite[Chapter II, \S7]{Borel}.

Another way to obtain special subdistributions/subfoliations is slope-stability. Given a coherent sheaf $F$ on $\pn$ of positive rank $\rk F$, the slope is defined by $\mu(F) :=c_1(F)/\rk F$. Then $F$ is called semistable if $\mu(F) \geq \mu(E)$ for every $E\subset F$ such that $\rk E <\rk F$. $F$ is called stable if the inequality is always strict. One also defines
\begin{align*}
    \mu_{\rm max}(F) &= \max \{\, \mu(E) \mid E \subset F \,\}; \\
    \mu_{\rm min}(F) &= \min \{\, \mu(Q) \mid F \twoheadrightarrow Q \, \text{and } Q \text{ is torsion free}  \,\}.
\end{align*}

We first record an observation used for $\fix(\sF)$. The condition $\mu_{\rm max}(F)\leq0$ holds, in particular, when $F$ is semistable with $c_1(F)=0$. 

\begin{lemma}\label{lem:sstab-h0}
Let $F$ be a rank-$r$ reflexive sheaf on $\pn$ such that $\mu_{\rm max}(F) \leq 0$. Then the evaluation morphism $H^0(F) \otimes \opn \hookrightarrow F$ is injective with saturated image. In particular, $h^0(F) \leq r$. If equality holds, then $F = \opn^{\oplus r}$. 
\end{lemma}

\begin{proof}
Let $h = h^0(F)$ and consider the evaluation map $\phi \colon \opn^{\oplus h} \to F$ generated by linearly independent global sections. We get a short exact sequence
\[
0 \lra \ker \phi \lra \opn^{\oplus h} \lra \im \phi \lra 0.
\]
Note that $\im \phi$ must be torsion-free. Hence, $\ker \phi \subset \opn^{\oplus h}$ is saturated. Since $\opn^{\oplus h}$ is semistable, $\mu(\ker \phi) \leq 0$ and $\mu(\im \phi) \geq 0$. On the other hand, $\mu(\im \phi) \leq \mu_{\rm max}(F) \leq 0$. Then $\mu(\im \phi) = \mu(\ker \phi) = 0$. By \cite[Corollary 1.6.11]{HL-sheaves}, $\ker \phi$ must be a direct summand of $\opn^{\oplus h}$. One can interpret the map $\ker \phi = \opn^{\oplus a} \hookrightarrow \opn^{\oplus h}$ as $\CC$-linear relations among the global sections defining $\phi$. Therefore, $\ker \phi = 0$, i.e., $\phi$ is injective, and $h\leq \rk F$. Consider $G$ the saturation of $\im \phi$ in $F$. Then we have an exact sequence
\[
0 \lra \opn^{\oplus h} \lra G \lra T \lra 0,
\]
where $T$ is a torsion sheaf whose first Chern class is the class of the divisorial part of its support. On the other hand, $\mu(G)\leq\mu_{\rm max}(F)\leq0$, and the exact sequence gives $c_1(T)=c_1(G)\leq0$. Hence $c_1(T)=0$, so $T$ is supported in codimension at least two. Since $G$ is reflexive, $G\cong \opn^{\oplus h}$ and $\im \phi = G$. If $h = \rk F$, then $G=F$ gives the last part of the statement.
\end{proof}

\begin{proposition}\label{prop:sstable-subfol}
Let $\sF$ be a foliation on $\pn$ such that $\mu_{\rm max}(\tF) \leq 0$. Then it admits a subfoliation of dimension $h^0(\tF)$ generated by global vector fields.
\end{proposition}

\begin{proof}
By Lemma \ref{lem:sstab-h0}, the evaluation map  $\phi \colon H^0(\tF) \otimes \opn \to \tF$ is injective with saturated image. Then $H^0(\tF)\otimes \opn \hookrightarrow \tF \hookrightarrow \tpn$ defines a subdistribution, which is a foliation since $H^0(\tF) = \fix(\sF)$ is closed under the Lie bracket. 
\end{proof}

Next, we provide a simple technical observation that will be useful to constrain the singular locus of a subfoliation $\sG$ whose tangent sheaf is locally free. Exact sequences, such as the one below, are known in Commutative Algebra as Bourbaki sequences; see, for instance, \cite{HKS} and the references therein. 

\begin{lemma}\label{lem:codim-ext}
Let $X$ be a smooth projective variety. Consider an exact sequence of sheaves on $X$
\[
0 \lra E \lra F \lra \sI_W(a) \lra 0
\]
where $E$ is locally free, $F$ is reflexive, and $\codim W \geq 2$. Then either 
\begin{itemize}
    \item $W= \emptyset$ (hence $F$ is locally free), or
    \item $W$ is pure of codimension two and Cohen-Macaulay off $\sing(F)$. 
\end{itemize}

\end{lemma}

\begin{proof}
Since $F$ is reflexive, it is locally free off $Z = \sing(F)$ whose codimension is at least three. Note that $Z\subset W$ (set-theoretically) since $F$ must be locally free at any point off $W$. At any point $p\in W\setminus Z$, the Hilbert-Burch theorem says that $W$ is Cohen-Macaulay of codimension two. 

We only need to prove that $W$ is pure or, equivalently, that $\OO_W$ satisfies the Serre condition $S_{1,2}$, see \cite[Proposition 1.1.10]{HL-sheaves}. This condition means that $\depth(\OO_{W,p}) \geq 1$ for any point $p\in W$ of codimension at least $3$. Fix such a point $p$. Since $F$ is reflexive and $E$ is locally free, we have $\depth(F_p) \geq 2$ and $\depth(E_p) =\codim p \geq 3$. Then $\depth(\sI_{W,p}) \geq 2$ and $\depth(\OO_{W,p}) \geq 1$, see \cite[Proposition 1.2.9]{BH-CMrings}.
\end{proof}

\begin{proposition}\label{prop:sing2xsplit}
Let $X$ be a smooth projective variety and $E$ a locally free sheaf on $X$. Let $F$ and $G$ be reflexive sheaves such that $G\subset F\subset E$, $F/G$ is torsion free of rank one, $G$ is locally free and  $H^1(G(c_1(G) - c_1(F))) = 0$. Then either 
\begin{itemize}
    \item $F \cong G \oplus \OO_X(c_1(F) - c_1(G))$, or
    \item $F/G = \sI_W(c_1(F) - c_1(G))$ with $W \subset \sing(E/G)$ a subscheme of pure codimension two.
\end{itemize} 
\end{proposition}

We will apply this proposition to the cases $E = \tpn$ and $E = \tpn(-\log D)$. 

\begin{proof}
By hypothesis, $F$ and $G$ fit in the following commutative diagram with exact rows and columns, where $\delta = c_1(F) - c_1(G)$.
\begin{equation}\label{eq:diagsubfol}
\begin{tikzcd}
 G \arrow[r, equal ] \arrow[d, hook, "\rho"] & G \arrow[d, hook, "\phi"] & \\
 F \arrow[r,hook, "\psi"] \arrow[d,two heads] & E \arrow[d,two heads] \arrow[r] & E/F \arrow[d,equal] \\
 \sI_W(\delta) \arrow[r,hook]  & E/G \arrow[r]  & E/F
\end{tikzcd}
\end{equation}
By Lemma \ref{lem:codim-ext}, either $W= \emptyset$ or $W$ is pure of codimension two. If $W = \emptyset$, then $F$ corresponds to an extension class in $\Ext^1(\OO_X(\delta) , G) = H^1(G(-\delta))$ which is zero by hypothesis. Thus $F = G \oplus \OO_X(\delta)$ in this case. 

Assume $W \neq \emptyset$. By Lemma \ref{lem:codim-ext}, $W$ is pure of codimension two. Let $p$ be an associated point of $\OO_W$. We have that $p$ has codimension two and, on the other hand, $F$ is reflexive. Thus $F_p$ is a free $\OO_{X,p}$-module. Of course, the same holds for $G$ and $E$. Then, we may regard $\phi_p = \psi_p\rho_p$ as a product of matrices. 

The ideal sheaf $\sI_Z$ of $Z= \sing (E/G)$ is defined by the $s\times s$ minors of $\phi$, where $s = \rk G$. Also note that $\sI_{W,p}$ is defined by the $s\times s$ minors of $\rho_p$. Then, by Binet-Cauchy Theorem, \cite[p.11]{Gantm.vol1}, each $s\times s$ minor of a product of matrices is a combination of $s\times s$ minors of both matrices. Hence, $\sI_{Z,p} \subset \sI_{W,p}$. Since $p$ is arbitrary, we get $\sI_Z \subset \sI_W$, see \cite[\S6]{Matsumura}. Therefore, $W\subset \sing (E/G)$.
\end{proof}

\begin{remark}\label{rem:integ-subfol}
1) Consider $\sF$ and $\sG$ as in the statement of Proposition \ref{prop:sing2xsplit}. Note that, in the first case, $\tF = \tG \oplus \opn(\delta)$ is a subsheaf of $\tpn$, hence $\delta \leq 1$, i.e., $c_1(\tF) \leq c_1(\tG)+1$. This is a simple constraint for having a subdistribution that splits. 

\noindent 2) If $\sF$ is integrable, then the integrability of $\sG$ is measured by the map $\extp^2\tG \to \tF/\tG$ that sends $v\wedge w$ to $[v,w] + \tG$. Since $\tG$ is locally free, this produces a section in $H^0(\extp^2\tG^\vee(\delta)\otimes \sI_W)$, with $W$ possibly empty. Then, $\sG$ is automatically integrable if this cohomology vanishes.  
\end{remark}

Combining Proposition \ref{prop:sstable-subfol} and Proposition \ref{prop:sing2xsplit}, one gets the following result.
\begin{corollary}\label{cor:sub-sing-cod2}
Let $\sF$ be a degree-$d$ foliation of dimension $k$ on $\pn$. If $h^0(\tF)=k-1$ and $\mu_{\rm max}(\tF)\leq0$, then either $\tF \cong \opn^{\oplus k-1} \oplus \opn(k-d)$ or $\fix(\sF)$ defines a subfoliation of dimension $k-1$ and with singularities in codimension two.
\end{corollary}

\subsection{2-dimensional foliations tangent to a global vector field} 
We now specialize to two-dimensional foliations such that $h^0(\tF) \geq  1$. If $h^0(\tF(-1)) > 0$ then $\sF$ is a linear pullback of a degree-$d$ one-dimensional foliation on $\p{n-1}$. Assuming $h^0(\tF(-1)) = 0$, a global section of $\tF$ defines a subfoliation of degree one. Then Proposition \ref{prop:sing2xsplit} applies. We start by proving Proposition~\ref{prop:B}, which generalizes \cite[Theorem C]{C} to foliations on $\pn$, $n\geq 3$.

\begin{proof}[Proof of Proposition~\ref{prop:B}]
If $h^0(\tF(-1))>0$ then $\sF$ is a linear pullback and $\tF \cong \opn(1)\oplus \opn(1-d)$. Suppose that $h^0(\tF(-1)) = 0$ and $\tF \not\cong \opn \oplus \opn(2-d)$. Since $\rk \tF = 2$, it is enough to prove that $c_3(\tF) \neq 0$. Consider the exact sequence from the first column of \eqref{eq:diagsubfol}:
\[
0  \lra  \opn  \lra  \tF  \lra  \sI_W(2-d)  \lra  0
\]
By Proposition \ref{prop:sing2xsplit} $W$ has pure codimension $2$. Furthermore, if $n>3$, then $W\cong \p{n-2}$ is a linear subspace; and if $n=3$, then $W$ is either one line or a (possibly degenerate) disjoint union of two lines, see \cite[\S2.3]{GJM} and Lemma \ref{lem:veccd2} below. 

If $n\geq 4$, we take a general $L \cong \p3$ linearly embedded in $\pn$ transverse to $W$. The above sequence then restricts to $L$ as
\[
0  \lra  \OO_L  \lra  \tF|_L  \lra  \sI_{W\cap L}(2-d)  \lra  0.
\]
$W\cap L$ is a line, and  $c_i(\tF|_L) = c_i(\tF)$ for $i = 1,2,3$. Therefore, we may assume that $n=3$. 

By \cite[Theorem 4.1]{H2}, $c_1(\tF) = 2-d$, $c_2(\tF) = \deg( W)$ and $c_3(\tF) = 2p_a(W)-2 + \deg(W)(d+2)$. If $W$ is a line then $\deg(W) = 1$ and $p_a(W) = 0$, hence $c_3(\tF) = d \geq 1$. If $W$ is (possibly the limit of) a disjoint union of two lines, $\deg(W) = 2$ and $p_a(W) = -1$, hence $c_3(\tF) = 2d\geq 2$. 
\end{proof}

We classify the possible subfoliations using Jordan normal forms.

\begin{lemma}\label{lem:veccd2}
Let $v\in H^0(\tpn)$ have zeros in codimension two. Then, up to a linear change of coordinates, either
\begin{enumerate}
 \item $v = x_0\del{ x_0} + \lambda x_1\del{ x_1}$ with $\lambda \in \mathbb{C}^*$; 
 \item $v = (x_0+ x_1)\del{ x_0} + x_1\del{ x_1}$;
 \item $v = x_1\del{ x_0} + x_3\del{ x_2}$;
 \item $v = x_1\del{ x_0} + x_2\del{ x_1}$;
 \item { $v = x_1\del{ x_0} + x_2\del{ x_2}$.} 
\end{enumerate}
\end{lemma}

\begin{lemma}\label{lem:sym-invar}
Let $\sF$ be a foliation on $\pn$. Let $v \in \aut(\sF) \setminus \fix(\sF)$ satisfy $v = \ell v'$ for some linear polynomial $\ell$. Then $H = V(\ell)$ is $\sF$-invariant.
\end{lemma}

\begin{proof}
Let $\omega$ be a $p$-form defining $\sF$, where $p = \codim \sF$. By hypothesis, $\lieder_v\omega = a\, \omega$ for some $a\in \CC$, and $v \cont \omega \neq 0$. Take $w \in H^0(\tF(t))$ for some $t\geq 0$. Then $\bkt{v,w} \cont \omega = 0$. Write $v \cont \omega = f\eta$, where $f$ is the greatest common divisor of the coefficients of $v\cont\omega$; note that $\ell\mid f$. Then $\lieder_w\omega  = g\, \omega$ for some polynomial $g$ and
\[
\lieder_w (v \cont \omega) = v\cont \lieder_w \omega  + \bkt{w,v} \cont \omega  = g v\cont \omega. 
\]
On the other hand, $w \cont \eta = 0$ hence 
\[
\lieder_w (v \cont \omega) = \lieder_w (f\eta) = w\cont d(f\eta) = w(f)\eta + f(w\cont d\eta).
\]
 Since the coefficients of $\eta$ have no common factor, we conclude that $f$ divides $w(f)$. Writing $f=\ell^m h$, with $\ell\nmid h$, we obtain $\ell\mid w(\ell)$. Since $w$ is arbitrary, $H$ is $\sF$-invariant.
\end{proof}

\begin{proposition}\label{prop:sub-deg2-alg}
Let $\sF$ be a foliation of dimension $2$ and degree $d\geq1$ on $\pn$. Assume that $h^0(\tF(-1))=0$, $h^0(\tF)>0$ and $\tF\not\cong\opn\oplus\opn(2-d)$. Then every nonzero $s\in H^0(\tF)$ induces an algebraic action singular in codimension 2. 
\end{proposition}

\begin{proof}
Let $\sG$ be the foliation induced by a nonzero section $s \in H^0( \tF)$. By Proposition \ref{prop:sing2xsplit}, we know that $\sG$ has singularities in codimension two. Suppose for contradiction that the general leaf of $\sG$ is not algebraic. Then, due to Lemma \ref{lem:veccd2}, $\sG$ is defined by either $v = x_0\del{ x_0} + \lambda x_1\del{ x_1}$ with $\lambda \not \in \QQ$, or $v = (x_0+ x_1)\del{ x_0} + x_1\del{ x_1}$, or $v = x_1\del{ x_0} + x_2\del{ x_2}$. Since $h^0(\tF(-1)) = 0$, we have the following.  

\begin{enumerate}
    \item If $v = x_0\del{ x_0} + \lambda x_1\del{ x_1}$ with $\lambda \not \in \mathbb{Q}$, then 
    \[
    x_0 \del{ x_0}, x_1 \del{ x_1} \in \afix(\sF) \setminus \fix(\sF);
    \]
    \item if $v = (x_0+ x_1)\del{ x_0} + x_1\del{ x_1}$, we have that 
    \[
    x_0\del{ x_0} + x_1\del{ x_1} ,x_1\del{ x_0} \in \afix(\sF)\setminus \fix(\sF); \quad \text{and}
    \]
    \item if $v = x_1\del{ x_0} + x_2\del{ x_2}$, we have that
    \[
    x_2\del{ x_2} ,x_1\del{ x_0} \in \afix(\sF)\setminus\fix(\sF).
    \]
\end{enumerate}
By Lemma \ref{lem:sym-invar}, $\sF$ has an invariant divisor $D$ composed of one or two coordinate hyperplanes. Hence $\tF \hookrightarrow \tpn$ factors through the logarithmic tangent sheaf $\tpn(-\log(D))$. By Proposition \ref{prop:sing2xsplit}, $v$ corresponds to a logarithmic vector field vanishing in codimension two. We will show this is absurd, concluding the proof. 

In the first case, $\tpn(-\log(D)) \cong \opn \oplus \opn(1)^{\oplus n-1}$ generated by $x_0\del{x_0} - x_1\del{x_1}$ and $\del{x_j}$, $j\geq 2$. Then the logarithmic vector field is 
\[
v - \frac{v(x_0x_1)}{2x_0x_1}\rad =  \frac{1-\lambda}{2} \left(x_0\del{x_0} - x_1\del{x_1} \right) - \frac{\lambda+1}{2}\sum_{j\geq 2} x_j\del{x_j},
\]
whose vanishing locus, as a section of $\tpn(-\log(D))$, is empty, since the coefficient $\frac{1-\lambda}{2}$ is not zero. The third case is analogous. For the second case, $\tpn(-\log(D)) \cong \opn(1)^{\oplus n}$ generated by $\del{x_j}$, $j \neq 1$. Then the logarithmic vector field is
\[
v - \rad = x_1 \del{x_0} - \sum_{j\geq 2} x_j\del{x_j},
\]
whose vanishing locus, as a section of $\tpn(-\log(D))$, is the point $(1:0:\cdots:0)$. 
\end{proof}

In light of the previous proposition, we need to understand foliations tangent to algebraic actions with codimension-two fixed points. We start with additive $\CC$-actions, which correspond to cases 3) and 4) in Lemma \ref{lem:veccd2} up to a linear change of coordinates.

\begin{proposition}\label{prop:class-om-nilp}
Let $\sF$ be a foliation of dimension $2$ and degree $2$ on $\pn$ tangent to a global $\CC$-action with codimension-two fixed points. Then $\sF \in \logg{1^{n-1},2}{2}{n}$.
\end{proposition}

\begin{proof}
We fix an action tangent to a vector field in either 3) or 4) in Lemma \ref{lem:veccd2}. Both cases follow the same argument, so we only prove one. In case 3), $v = x_1\del{ x_0} + x_3\del{ x_2}$. Then it has a rational first integral $\rho\colon \pn \dashrightarrow \PP(1^{n-1},2)$ given by $\rho(x_0: \dots : x_n) = (x_1:x_3:x_4:\dots :x_n: x_0x_3-x_1x_2)$. Then $\sF = \rho^*\sH$ for $\sH$ a foliation by curves on $\PP(1^{n-1},2)$. Note that $\rho$ does not ramify in codimension one, hence 
\[
\opn(n+1) = \det(\nF) = \rho^*\det(\nH) = \opn(n+\deg \sH).
\]
Thus $\deg(\sH) = 1$ and $\sH$ is defined by a global vector field. As we observed in Remark \ref{rem:log(1,2)}, $\sF \in \logg{1^{n-1},2}{2}{n}$. 
\end{proof}

Actions of the multiplicative group $\CC^\ast$ with a codimension-two fixed locus are,  up to a linear change of coordinates, induced by vector fields as in 1) in Lemma \ref{lem:veccd2}, $v = x_0\del{ x_0} + \mu x_1\del{ x_1}\in H^0(\tpn)$ with $\mu \in \QQ$. The next result classifies the $2$-dimensional foliations invariant under such actions. Note that we don't impose any restriction on $\mu\in \CC^*$ a priori. 

\begin{proposition}\label{prop:class-om-diag}
Let $\sF$ be a two-dimensional foliation on $\pn$ of degree $2$. Suppose that $v = x_0\del{ x_0} + \mu x_1\del{ x_1}\in H^0(\tF)$ with $\mu \neq 0$. Then either $\sF$ is contained in a logarithmic component or $\sF$ is induced by an action of $\mf{aff}(\CC)$.
\end{proposition}

\begin{proof}
Let $\omega \in H^0(\Omega_{\pn}^{n-2}(n+1))$ be a form defining $\sF$. Then $v\cont \omega = 0$ and $\lieder_v\omega = a\omega$ for some $a\in\CC$. Let $v'= (n+1)v - a\rad$. Then $v'\cont \omega = v'\cont d\omega = 0$. 

If $a\neq 0$ then, the coefficients of $v'$ vanish in codimension at least $3$ in $\CC^{n+1}$, unless $n =3$,  $a = 4$ , and $\mu = 1$. In this case, the weight $a=4$ implies that $\omega = q(x_0,x_1)(x_0dx_1-x_1dx_0)$ for some $q$ of degree $2$. This is not possible since $\omega$ cannot vanish in codimension one. 

By the dual Saito Division Lemma, see \cite[Corollary 4.2.]{C}, there exists $\eta$ such that $d\omega = v'\cont \eta$. Since $\eta$ is an $n$-form in $n+1$ variables, there exists a linear vector field $u$ such that $\eta = u \cont \dV$, where $\dV = dx_0\wedge \dots \wedge dx_n$. Then
\[
\omega = \rad \wedge v \wedge u \cont \dV. 
\]
Hence $\tF \cong \opn^{\oplus 2}$ is globally generated by $v$ and $u$. If $u$ and $v$ commute, then $\sF \in \logg{1^{n+1}}{2}{n}$, as in Remark \ref{rem:log(1,2)}. If they don't commute, then $\sF$ is defined by an action of $\mf{aff}(\CC)$. By Remark \ref{Alg: log(1,2)}, $\sF$ is either in an affine-Lie component or in a logarithmic component. 

Now suppose that $a = 0$, so that $v\cont \omega = v\cont d\omega  = 0$. Then 
\[
d\omega =  (x_0 dx_1 - \mu x_1 dx_0) \wedge \alpha + \beta, 
\]
where $\alpha$ and $\beta$ do not depend on $dx_0$ or $dx_1$. Hence, $\beta$ is an $(n-1)$-form in $n-1$ variables, i.e., $\beta = q \, \Theta$ for $q$ a quadratic polynomial and $\Theta=dx_2\wedge\dots\wedge dx_n$. Moreover, $\alpha = u \cont \Theta$ for some linear vector field $u$. Hence,
\[
0 = \lieder_v d\omega = (x_0 dx_1 - \mu x_1 dx_0) \wedge ((\bkt{v,u} + (\mu+1)u) \cont \Theta) + v(q)\, \Theta.
\]
Thus, $v(q) = 0$ and $\bkt{v,u} + (\mu+1)u = 0$. Note that $u = \sum_{j\geq 2}b_j \del{x_j}$ for some linear polynomials $b_j$. Hence, we have $v(b_j) = -(\mu+1)b_j$ for every $j\geq 2$. This forces one of the following cases: $\mu = -1$ and $\frac{\partial b_j}{\partial x_0} = \frac{\partial b_j}{\partial x_1} = 0$; $\mu = -2$ and $u = x_0u'$ for $u'$ constant; or $\mu = -\frac{1}{2}$ and $u = x_1u'$ for $u'$ constant. The last two cases are identified up to exchanging $x_0$ and $x_1$.

If $\mu = -1$, then $q = c x_0x_1 + q_0(x_2, \dots, x_n)$ for some $c\in\CC$ and a quadratic polynomial $q_0$. Then 
\[
 (n+1)\omega = \rad \cont d\omega  =  -d(x_0x_1) \wedge (\rad \wedge u \cont  \Theta) + 2 x_0x_1 (u \cont \Theta) \, + \,q \, \rad \cont \Theta.
\]
It follows that $\sF$ is the pullback of a foliation on $\PP(2,1^{(n-1)})$ via the rational map $\pn \dashrightarrow \PP(2,1^{(n-1)})$ given by $y_1 = x_0x_1$, $y_j = x_j$ for $j\geq 2$. By Remark \ref{rem:log(1,2)}, $\sF \in \logg{1^{n-1},2}{2}{n}$.

If $\mu = -2$, then $q =  q(x_2, \dots, x_n)$ and, up to a linear change of coordinates, $u = x_0\del{x_2}$. Then
\[
d\omega =  d(x_0^2x_1) \wedge dx_3\wedge \dots \wedge dx_n + q\, \Theta = ( d(x_0^2x_1) + qdx_2) \wedge dx_3\wedge \dots \wedge dx_n.
\]
Then by the Poincaré Lemma, there exists a polynomial $f$ such that $d\omega = df \wedge dx_3\wedge \dots \wedge dx_n$. Therefore $\sF \in \logg{1^{n-2},3}{2}{n}$.
\end{proof}

\section{Stable two-dimensional foliations}\label{sec:stable-dim2} In this section we describe two-dimensional foliations whose tangent sheaf is stable. We prove that the foliation is either rational or of a special type that we call a \emph{null-correlation} foliation. The main ingredient is the cokernel of  $\tF \hookrightarrow \opn(1)^{\oplus n+1}$, given by Lemma \ref{lem:lift-map}.

As defined in \cite{BeRo}, the \emph{global-sections type} of a torsion-free sheaf $F$ of rank $r$ is the set of integers $(a_1, \dots, a_r)$ such that $a_1 \geq \dots \geq a_r$ and there exists an injective morphism $\phi \colon \bigoplus_{j=1}^r\opn(a_j) \hookrightarrow F$ that is maximal in the following sense. If $\psi \colon \bigoplus_{j=1}^s\opn(b_j) \hookrightarrow F$ is an injective map for some $b_1\geq \dots \geq b_s$, then $s\leq r$ and $b_j \leq a_j$ for every $j$. We write $\gst(F) = (a_1, \dots, a_r)$ for its global-sections type. If $F$ is generically globally generated, i.e., there is a morphism $\opn^{\oplus k} \to F$ surjective on some open subset, then $a_r \geq 0$ \cite[Lemma 6]{BeRo} and there are bounds for the Chern classes of a reflexive $F$, see \cite[Theorem 1]{BeRo}:
\begin{equation}\label{eq:BeRo-thm1}
    c_1(F) \geq \sum_j a_j , \quad c_2(F) \geq \sum_{i<j} a_ia_j,
\end{equation}
and equality in either bound yields $F \cong \bigoplus_{j=1}^r\opn(a_j)$. They also prove a more general bound, see \cite[Theorem 2]{BeRo}:
\begin{equation} \label{eq:BeRo-thm2}
    c_2(F) \geq \sum_{i<j} a_ia_j + \left(c_1(F) - \sum_{j} a_j\right)\left( \sum_j a_j +1-a_2\right).
\end{equation}

Let $\sF$ be a foliation of degree two and dimension $k\geq 2$ on $\pn$, defined by $\omega\in H^0(\Omega_{\pn}^{n-k}(n-k+3))$. The contraction against $d\omega$ defines a map $\psi \colon \opn(1)^{\oplus n+1}= \widetilde{T_S} \to \widetilde{\Omega_S^{n-k}}(n-k+3)$ whose image we call $M$. By Lemma \ref{lem:lift-map}, we obtain an exact sequence 
\begin{equation}\label{eq:stable-main-seq}
0\lra\tF\lra\opn(1)^{\oplus n+1}
\stackrel{\psi}{\lra}M^{\vee\vee}\lra Q\lra0,
\end{equation}
where $Q=M^{\vee\vee}/M$ is supported in codimension at least two.

\begin{lemma}\label{lem:reduc-stable-case}
Let $\sF$ be a degree-two foliation satisfying $H^0(\tF(-1)) = H^0(\tF^\vee)=0$. Then either
\begin{enumerate}
\item $M^{\vee\vee}\cong\opn(3)\oplus\opn(1)^{\oplus n-k}$ or $ M^{\vee\vee}\cong\opn(2)^{\oplus2}\oplus\opn(1)^{\oplus n-k-1}$; or
\item $k=2,3$ and $\tF=\pi^*F$, where $\pi\colon \pn\dashrightarrow\p3$ is a linear projection and $F$ is either $\tp3(-1)$ or a null-correlation bundle $N$ on $\p3$. In this case, the non-Kupka locus has codimension at least three.
\end{enumerate}
\end{lemma}

Since $\pi$ is rational, $\pi^*N$ means the reflexive extension $j_*(\pi|_U^*N)$, where $U$ is the domain of definition of $\pi$ and $j\colon U\hookrightarrow\pn$ is the inclusion. Note that $\codim(\pn \setminus U) = 4$. Note that the cohomological hypotheses of the lemma hold if $\tF$ is semistable and $k\geq 3$, or stable and $k=2$. 

\begin{proof}
Note that the support of $Q$ has codimension at least two and agrees, in codimension two, with the vanishing locus of $d\omega$, i.e., the non-Kupka locus. We set $\delta=\ch_2(Q)\geq0$.  Then $\delta=0$ if and only if $\codim\Supp Q\geq3$. Let $F=M^{\vee\vee}(-1)$. From \eqref{eq:stable-main-seq} and $c_1(\tF)=k-2$ we obtain $c_1(F)  = 2$ and, since $\tF$ is locally free in codimension two, 
\[
c_2(\tF^\vee(1)) = c_2(\tF(-1)) = 4 - c_2(F)  -\delta.
\]
Since $\psi$ is generically surjective, $F$ is generically globally generated. The global-sections type $\gst(F)=(a_1,\dots,a_{n+1-k})$ satisfies $a_1 \geq \dots \geq a_{n+1-k}\geq 0$  and $\sum_j a_j \leq 2$, due to \eqref{eq:BeRo-thm1}, then $a_3=\cdots=a_{n+1-k}=0$ and 
\[
(a_1,a_2)\in\{(2,0),(1,1),(1,0),(0,0)\}.
\]
If $(a_1,a_2)=(2,0)$ or $(1,1)$, then $F$ splits, by \cite[Theorem 1]{BeRo}, giving the first item. 

Assume that $F$ does not split. Then $(a_1,a_2)$ is $(1,0)$ or $(0,0)$, and \eqref{eq:BeRo-thm2} yields $c_2(F)\geq2$. Thus $c_2(\tF^\vee(1))\leq 2 -\delta$. On the other hand, $\tF^\vee(1)$ is also nonsplit, generically globally generated with $c_1(\tF^\vee(1)) = 2$. The same argument shows that $c_2(\tF^\vee(1)) \geq 2$. Therefore, $\delta = 0$. 

From $\delta = 0$, $\tF^\vee(1)$ is globally generated off a set of codimension at least $3$. By \cite[Lemma 2]{BeRo}, $k-1$ general global sections give
\[
0\lra\opn^{\oplus k-1}\lra\tF^\vee(1)\lra\sI_Y(2)\lra0,
\]
where $Y$ is of pure codimension $2$ and degree $2$, and $Y$ is smooth in codimension $2$. Note that $h^0(\sI_Y(1)) = h^0(\tF^\vee) = 0$, by hypothesis. Hence, $Y$ is nondegenerate. This forces $Y$ to be reducible, i.e., $Y = L_1 \cup L_2$, for $L_i$ a linear subspace. Up to a linear change of coordinates, $Y = V((x_0,x_1)\cap (x_2,x_3))$. This yields the minimal resolution 
\begin{equation}\label{eq:res-nullcorr}
0\lra\opn(-3)\stackrel{P_2}{\lra}\opn(-2)^{\oplus4}
\stackrel{P_1}{\lra}\opn(-1)^{\oplus k+3}\lra\tF^\vee\lra0.
\end{equation}
The entries of $P_2$ are four independent linear polynomials. Up to a change of coordinates, we may take $P_2=[x_0\, x_1\, x_2 \,x_3]^T$. From $P_1P_2 = 0$, the entries of $P_1$ must also depend only on the first 4 variables, and this also forces $k\leq 3$; otherwise, $\tF^\vee(1)$ would have a trivial summand, but $h^0(\tF(-1)) = 0$ by hypothesis. Hence, $\tF^\vee$ is the reflexive pullback of a rank-$k$ sheaf on $\p3$ with the same resolution. If $k = 3$, this is $\Omega_{\p3}^1(1)$, by Koszul. If $k=2$, by \cite[Lemma 4.3.2]{OkonekSchneiderSpindler1980}, this sheaf is a null-correlation bundle. 
\end{proof}

Fix an integer $0 \leq m \leq n$. Then the exterior derivative $d\colon S \to \Omega_S^1$ can be written as $d = d' + d''$ where $d'$ differentiates with respect to the first $m+1$ variables and $d''$ acts on the others. For $R = \CC[x_{m+1}, \dots, x_n]$, $d'$ can be seen as the relative differential of $S$ over $R$. Then for any form $\alpha$ we write $\alpha = \alpha' + \alpha''$ where $\alpha'$ is the summand comprising only the differentials $dx_i$ for $i\leq m$. We have that 
\[
d\alpha \wedge dx_{m+1}\wedge \dots \wedge dx_n = d'\alpha' \wedge dx_{m+1}\wedge \dots \wedge dx_n.
\]
Moreover, since $d'\alpha'$ only involves $dx_i$ with $i\leq m$, the form above vanishes if and only if $d'\alpha' = 0$. Furthermore, there is a relative Poincaré Lemma: 
\begin{equation}\label{eq:poincare}
    d'\alpha' = 0 \Longrightarrow \exists \beta \text{ such that } \alpha' = d'\beta.
\end{equation}

\begin{lemma}\label{lem:split-to-rational}
Let $\sF$ be a degree-two foliation of dimension $k\geq 2$ such that $H^0(\tF^\vee) = 0$. With the notation as above, if $M^{\vee\vee}$ splits, then $\sF$ is a rational foliation.
\end{lemma}

This lemma holds, in particular, for semistable $\tF$ if $k>2$ or stable $\tF$ if $k=2$. 

\begin{proof}
By Lemma \ref{lem:decomp-free},  $\sF$ is defined by a form $\omega$ such that $d\omega=\theta_1\wedge\cdots\wedge\theta_{n+1-k}$, where each $\theta_i$ has total degree $\deg \theta_i \geq 1$ and $\sum \deg \theta_j = n+3-k$. Then  $\theta_j$ has total degree $1$ for $j\geq 3$ and the total degrees of $\theta_1,\theta_2$ are $(3,1)$ or $(2,2)$. Up to a linear change of coordinates, we assume that $\theta_j = dx_{j+k-1}$ for $j\geq 3$. Moreover, $\theta_1$ and $\theta_2$ only involve the first $k+2$ differentials. 

Suppose first that the total degrees are $(3,1)$. We may take $\theta_2=dx_{k+1}$ and
\[
d\omega=\theta_1\wedge dx_{k+1}\wedge\cdots\wedge dx_n.
\]
Let $d'$ be the derivative with respect to the first $k+1$ variables. Then $d'\theta_1 = 0$ and, by \eqref{eq:poincare}, there exists a homogeneous polynomial $f$ of degree $3$ such that $d'f = \theta_1$. Hence, $\sF \in \logg{1^{n-k},3}{k}{n}$. 

Assume now that both $\theta_1$ and $\theta_2$ are of degree two. Let $d'$ be the derivative with respect to the first $k+2$ variables. Then
\[
d\omega=\theta_1\wedge\theta_2\wedge dx_{k+2}\wedge\cdots\wedge dx_n,
\]
with $\theta_1,\theta_2$ involving only $dx_0,\dots,dx_{k+1}$, and $d'(\theta_1\wedge\theta_2)=0$.
It is enough to show that $d'\theta_1 = d'\theta_2 = 0$. In that case, $\theta_i = d'q_i$ for quadratic polynomials $q_1$ and $q_2$ hence   $\sF\in\logg{1^{n-k-1},2,2}{k}{n}$.

Suppose $d'\theta_1 \neq 0$. Then $d'\theta_1 \wedge \theta_1 \wedge \theta_2 = 0$. Consider 
\[
\gamma = \rad \cont (d\theta_1 \wedge \theta_2 \wedge dx_{k+2}\wedge \dots \wedge dx_n), \quad \rho = \rad \cont (d\theta_1 \wedge \theta_1 \wedge dx_{k+2}\wedge \dots \wedge dx_n). 
\]
These would define subdistributions of $\sF$ of dimension $k-1$ and degree at most one. But these are excluded by $H^0(\tF^\vee) = 0$. Hence, $\gamma = \rho = 0$. Then $d'\theta_1 \wedge \theta_1 =d'\theta_1 \wedge \theta_2 = 0$. Since $d'\theta_1$ has constant coefficients, it is proportional to $\theta_{1,p}\wedge\theta_{2,p}$ for a general closed point $p$. Therefore, regarding $\theta_{1,p}$ and $\theta_{2,p}$ as constant 1-forms, $\theta_1\wedge \theta_2 = g \,\theta_{1,p}\wedge\theta_{2,p}$ for some quadratic polynomial $g$. This means that $d\omega$ has zeros in codimension one, which is absurd. We conclude that $d'\theta_1 = 0$. Similarly, $d'\theta_2 = 0$.
\end{proof}

\subsection{Null-correlation foliations}\label{subsec:nullcorr-common} Next we describe the two-dimensional foliations in the second item of Lemma~\ref{lem:reduc-stable-case}. Assume that the projection $\pi$ is defined by the first 4 variables. Then we can assume $\tF$ has a resolution as in \eqref{eq:res-nullcorr} with $P_2=[x_0\, x_1\, x_2 \,x_3]^T$ and $k=2$. Up to a linear change of coordinates, we can assume that 
\begin{equation}\label{eq:matB}
    P_1 = \begin{bmatrix}
    -x_{1}&x_{0}&-x_{3}&x_{2}\\
      0&-x_{2}&x_{1}&0\\
      -x_{2}&0&x_{0}&0\\
      0&-x_{3}&0&x_{1}\\
      -x_{3}&0&0&x_{0}
\end{bmatrix}.
\end{equation}

\begin{lemma}\label{lem:reduction-nullcorr}
Let $\sF$ be a two-dimensional foliation as in the second alternative of Lemma~\ref{lem:reduc-stable-case}. Then the inclusion $\tF\hookrightarrow\opn(1)^{\oplus n+1}$ defines a subspace  $W\subset H^0(N(1))$. If $M^{\vee\vee}$ does not split, then either
\begin{enumerate}
\item $\dim W=5$ and $M^{\vee\vee}\cong\opn(1)^{\oplus n-4}\oplus\pi^*\Omega^1_{\p3}(3)$; or
\item $\dim W=4$ and $M^{\vee\vee}\cong\opn(1)^{\oplus n-3}\oplus\pi^*N'(2)$ for another null-correlation bundle $N'$.
\end{enumerate}
In either case, $Q$ is supported on the indeterminacy locus of $\pi$, which is linear of codimension 4.
\end{lemma}

\begin{proof}
By the resolution \eqref{eq:res-nullcorr}, the inclusion $\tF\hookrightarrow\opn(1)^{\oplus n+1}$ lifts to a map $P_0\colon \opn(-1)^{\oplus 5} \to \opn(1)^{\oplus n+1}$. Assume that the projection $\pi$ is defined by the first 4 variables. For $P_1$ as in \eqref{eq:matB}, the equation $P_0P_1 = 0$ implies that $P_0$ depends only on $x_0,\dots,x_3$. Hence, $M^{\vee\vee}=\pi^*G$ for a reflexive sheaf $G$ on $\p3$. We have
\[
0\lra N\lra \op3(1)^{\oplus n+1}\lra G\lra Q_0\lra0,
\]
with $\codim \Supp Q_0\geq3$. The leftmost map defines $n+1$ elements in $\Hom(N,\op3(1)) \cong H^0(N(1))$. Let $W \subset H^0(N(1))$ be the subspace spanned by these sections. If $r=\dim W$, then $G\cong\op3(1)^{\oplus n+1-r}\oplus G'$, where $\rk G'=r-2$. Since $G$ does not split, $r\geq4$, while $h^0(N(1))=5$, so $r\in\{4,5\}$. By construction, $N = \coker P_1$ for $P_1$ as in \eqref{eq:matB}. Then a basis of $H^0(N(1))$ is represented by the rows of
\begin{equation}\label{eq:nullcorr-B}
B=\begin{bmatrix}
0&-x_0x_3&x_1x_3&x_0x_2&-x_1x_2\\
x_0x_3&0&x_3^2&x_0^2&-x_0x_1-x_2x_3\\
-x_1x_3&-x_3^2&0&-x_0x_1+x_2x_3&x_1^2\\
-x_0x_2&-x_0^2&x_0x_1-x_2x_3&0&x_2^2\\
x_1x_2&x_0x_1+x_2x_3&-x_1^2&-x_2^2&0
\end{bmatrix}.
\end{equation}
This skew-symmetric matrix has constant rank two. If $r=5$, $G' \cong \coker B$. Dualizing
\[
0\lra N\to\op3(1)^{\oplus5}\lra G'\lra 0
\]
gives $G'\cong\Omega^1_{\p3}(3)$. If $r=4$, then we have the following commutative diagram
\[
\begin{tikzcd}
    & \op3(1) \ar[r,equal] \ar[d,hook] & \op3(1) \ar[d,hook] \\
    N \ar[d,equal]\ar[r,hook] & \op3(1)^{\oplus 5} \ar[d, two heads] \ar[r, two heads] & \Omega_{\p3}^1(3) \ar[d, two heads]\\
    N \ar[r,hook] & \op3(1)^{\oplus 4} \ar[r, two heads] & G_0
\end{tikzcd}
\]
where $G_0^{\vee\vee}  =G'$. Dualizing the rightmost column gives a degree-zero distribution on $\p3$ with a non-split tangent sheaf. By \cite[Proposition 7.1]{MOJ}, $G'\cong N'(2)$ for a null-correlation bundle $N'$.
\end{proof}

\begin{lemma}\label{lem:nullcorr-n>=7}
A null-correlation foliation satisfies $n\geq7$ and is tangent to the projection $\pi\colon \pn\dashrightarrow\p3$.
\end{lemma}

\begin{proof}
Let $\widehat\sF$ be the foliation on $\CC^{n+1}$ defined by $d\omega$, where $\omega$ defines $\sF$. The inclusion $\tF\hookrightarrow\opn(1)^{\oplus n+1}$  lifts to $P_0 \colon \opn(-1)^{\oplus 5} \to \opn(1)^{\oplus n+1}$. By \eqref{eq:res-nullcorr} and \eqref{eq:matB}, we have a factorization $P_0 = C^T B$, where $B$ is \eqref{eq:nullcorr-B} and $C=(c_{ij})$ is a constant $5\times(n+1)$ matrix. Put
\[
D_i=\sum_{j=0}^nc_{ij}\del{x_j},\quad
v_i=\sum_{j=0}^4b_{ij}D_j, \quad \text{and} \quad \Pi_C=\sum_{0\leq i<j\leq4}b_{ij}D_i\wedge D_j,
\]
where $B = (b_{ij})$. The rank-4 Pfaffians of $B$ vanish, hence we have the relations $b_{ij}b_{kl} - b_{ik}b_{jl}+b_{il}b_{jk} =0$, which give $v_i\wedge v_j=b_{ij}\Pi_C$. Then $\Pi_C$ defines $\widehat\sF$. The integrability is then read from the bracket $\bkt{v_i,v_j} \wedge v_i \wedge v_j = 0$ for any pair of indices $i,j$. Let us take $v_1$ and $v_2$. We have that 
\begin{equation}\label{eq:aux-v12}
    \bkt{v_1,v_2} = fv_1 + g v_2,
\end{equation}
for some rational functions $f$ and $g$. Taking the wedge against $v_1$ and $v_2$ yields that $fx_3^2$ and $gx_3^2$ are polynomials.

Since $\rk C \geq 4$, let $(k_0, \dots, k_4)$ span the (maybe trivial) kernel of $C^T$. Write $\bkt{v_1,v_2}  = \sum_j p_j D_j$. Then \eqref{eq:aux-v12} yields, for some rational function $h$, the equations
\begin{align*}
     k_0 h &= p_0 - fx_0x_3 +gx_1x_3, &\quad  k_1 h &= p_1 +gx_3^2, &\quad
      k_2 h &= p_2 -fx_3^2 ,\\
      k_3 h &= p_3 -fx_0^2 +g(x_0x_1-x_2x_3), &\quad 
     k_4 h &= p_4 +f(x_0x_1+x_2x_3) -gx_1^2.&&
\end{align*}
Since $fx_3^2$ and $gx_3^2$ are polynomials, $h$ is also a polynomial. We may write $f = (p_2 -k_2h)/x_3^2$ and $g = (-p_1+k_1h)/x_3^2$. Hence we get the following equations:
\begin{align*}
   (k_0x_3^2 -k_2x_0x_3 -k_1x_1x_3) h &= p_0x_3^2 - p_2 x_0x_3 -p_1x_1x_3 ,\\
   (k_3 x_3^2  -k_2x_0^2 -k_1(x_0x_1-x_2x_3)) h &= p_3x_3^2 -p_2 x_0^2 -p_1(x_0x_1-x_2x_3), \\
   (k_4x_3^2 +k_2(x_0x_1+x_2x_3)+k_1x_1^2) h &= p_4x_3^2 +p_2 (x_0x_1+x_2x_3) +p_1x_1^2.
\end{align*}
Note that $p_j = \sum_{i=0}^4 (b_{1i}D_i(b_{2j}) - b_{2i}D_i(b_{1j}))$. Substituting in the previous equations and collecting coefficients yields a linear system involving $D_i(x_j)$ for $0\leq j \leq 3$. We conclude that either
\[
D_i(x_j) = 0 \quad \text{for} \quad 0\leq j \leq 3, \, 0\leq i \leq 4,
\]
or $(k_0, \dots, k_4)\neq 0$ and satisfies $k_0^2 +k_1k_4-k_2k_3 =0$. By reducing to the cases $k_j \neq 0$, $j=1, \dots 4$, one concludes that in this case $\Pi_C$ degenerates along a linear subspace of codimension two, which contradicts Lemma \ref{lem:reduction-nullcorr}. Therefore, each $D_j$ depends only on the derivations $\del{x_i}$, $i\geq 4$. Equivalently, the first four columns of $C$ are zero. Since $\rk C \geq 4$ we have $n\ge 7$. Moreover, $\Pi_C$ is tangent to $\pi$. This concludes the proof.
\end{proof}

\begin{remark}\label{rem:nullcorr-geom}
Lemma \ref{lem:nullcorr-n>=7} also gives a geometric description of a null-correlation foliation $\sF$. Let $\phi \colon \pn \dashrightarrow \p3$ be the associated linear projection and $W$ its base locus. Each fiber of $\phi$ is a $\p{n-3} \supset W$ and the restriction of $\sF$ to the fiber is given by a constant bivector away from $W$. Hence, the leaves are open in linear embedded planes. 
\end{remark}

\section{The null-correlation component}\label{sec:nullcorr-component}
Lemma \ref{lem:reduction-nullcorr} divides the null-correlation foliations into two classes. We call those that fit into item (1) \emph{general null-correlation foliations} and those that fit into item (2) \emph{special null-correlation foliations}. The following lemma justifies this terminology. 

\begin{lemma}\label{lem:special-to-general}
Let $n\geq8$. Every special null-correlation foliation is a limit of general null-correlation foliations.
\end{lemma}

\begin{proof}
Up to a linear change of coordinates, we may assume we are as in the proof of Lemma \ref{lem:nullcorr-n>=7}. We only need to set the vector fields $D_j$. Up to a further change of coordinates, suppose that $D_0 = 0$ and $D_j = \del{x_{4+j}}$ for $j\geq 1$. Then define the deformation $\sF_t$ by setting $D_0(t) = t\del{x_{4}}$. For every $t\neq 0$ we have that the $D_j$ span a 5-dimensional space, hence $\sF_t$ is a general null-correlation foliation. 
\end{proof}

From the construction, general null-correlation foliations are projectively equivalent to each other. Indeed, we can choose $\pi$ to be the projection on the first 4 variables and change the last $n-3$ variables so that $D_i = \del{x_{4+i}}$. We then fix a foliation $\sF$ given by 
\begin{equation}\label{eq:Pi-general-nullcorr}
\Pi=\sum_{0\leq i<j\leq4}b_{ij}(x_0,\dots,x_3)
\del{x_{4+i}}\wedge\del{x_{4+j}},
\end{equation}
where $B=(b_{ij})$ is \eqref{eq:nullcorr-B}, and consider the action of $\SL_{n+1}$. 

\begin{lemma}\label{lem:orbit-nullcorr}
For $n\geq 8$,  the orbit of $[\Pi] \in \PP H^0(\extp^2\tpn)$ under the action of $\SL_{n+1}$ has dimension $9n-23$.
\end{lemma}
\begin{proof}
Let $V = \CC^{n+1}$. Note that the curl of $\Pi$ is zero, hence $\phi \in \SL(V)$ preserves $\sF$ if and only if $\phi_* \Pi = c\,\Pi$ for some $c\in \CC^*$. Since $\Pi$ depends only on $x_0, \dots, x_3$, and $\del{x_4}, \dots, \del{x_8}$, we may decompose $V=V_1\oplus V_2\oplus V_3$, with $V_1 = \langle e_0, \dots, e_3\rangle$, $V_2 = \langle e_4 , \dots, e_8 \rangle$, and $V_3 = \langle e_9 , \dots, e_n\rangle$, or $V_3 = 0$ if $n=8$. Let $G < \SL(V)$ be the stabilizer of $\Pi$. For $\phi \in G$, we have the block decomposition
\[
\phi = \begin{bmatrix}
    \phi_{00} & 0 & 0 \\
    \phi_{10} & \phi_{11} & \phi_{12} \\
    \phi_{20} & 0 & \phi_{22} 
\end{bmatrix}.
\]
Then $\phi_*\Pi = \Pi$ translates to $B ( \phi_{00}x) =  \phi_{11}B(x)\phi_{11}^T$. Linearizing the action at the identity gives the map
\[
(\alpha,\beta)\longmapsto
\frac{dB}{dx}(\alpha x)-\beta B-B\beta^T,
\qquad
\alpha\in \gl_4,\ \beta\in\gl_5.
\]
One then verifies directly (maybe in Macaulay2) that this linear system has rank $30$. The off-diagonal blocks are free, and $\det(\phi_{22})\det(\phi_{00})\det(\phi_{11}) = 1$. Hence, $\dim G = n^2 - 7n +22$. Therefore, the dimension of the orbit of the class $[\Pi] \in \PP H^0(\extp^2\tpn)$ is
\[
\dim(\SL(V)\cdot [\Pi]) = \dim(\SL(V)\cdot \Pi) - 1 = 9n -23. 
\]
\end{proof}

\begin{proposition}\label{prop:nullcorr-n>=8}
For $n\geq 8$, let $\NC_n$ be the closure of the orbit of $[\Pi]$ in $\Fol(2,2,n)$. Then $\NC_n$ is a generically reduced irreducible component of the space of foliations $\Fol(2,2,n)$. 
\end{proposition}

\begin{proof}
First, we identify $H^0(\extp^2\tpn)$ with the space of curl-free bivectors with quadratic coefficients. Consider a first-order deformation $\Pi_\epsilon = \Pi + \epsilon P \mod \epsilon^2$, where $P$ is a curl-free bivector with quadratic coefficients. By Lemma \ref{lem:curl-free-integ}, integrability reads $\Pi \wedge P = 0$. Thus,
\[
T_{[\Pi]}\Fol(2,2,n) \cong \{\, P \mid D(P) = \Pi\wedge P = 0 \,\}/ \CC\cdot\Pi.
\]
We will show that  $\Pi_\epsilon = (\id + \phi)_{\epsilon}\Pi$, hence $T_{[\Pi]}\Fol(2,2,n) = T_{[\Pi]}\NC_n$. 

Let $A = \CC[\epsilon]/(\epsilon^2)$ be the ring of dual numbers and $S_A = S\otimes A = A[x_0, \dots, x_n]$ be the extended polynomial ring. Let $M_0 = \{\,v\in T_S \mid v\wedge \Pi = 0\,\}$, so that $\widetilde{M_0} = \tF$. Consider the $S_A$-module $M_{\epsilon} = \{\, v \in T_S \otimes A \mid v\wedge \Pi_\epsilon = 0\,\}$. We will show that $M_\epsilon$ is a (flat) first-order deformation of $M_0$.

Any $v \in T_S \otimes A$ can be written as $v = v' + \epsilon v''$. Then $v \in M_\epsilon$ implies that $v'\in M_0$. Hence, we get a morphism $\pi \colon M_\epsilon \to M_0$ which sends $v$ to $v'$. Note that $M_0$ is generated by the vector fields $v_i = \sum_j b_{ij} \del{x_{j+4}} = -dx_{i+4}\cont \Pi$. These lift to $v_i^\epsilon = -dx_{i+4}\cont \Pi_{\epsilon}$. Thus, $\pi$ is surjective. Moreover, for $w \in T_S$, $\epsilon w \in M_{\epsilon}$ if and only if $w \in M_0$. By \cite[Corollary 6.2]{Eisen}, $M_\epsilon$ is flat over $A$. 

The $v_i^\epsilon$ generate $M_{\epsilon}$ by Nakayama's Lemma. Then we get a morphism $\psi_\epsilon\colon S_A(-1)^{\oplus 5} \twoheadrightarrow M_{\epsilon}$. By flatness, $K = \ker \psi_\epsilon$ is a flat deformation of $\ker \psi_0$. Iterating, we get a resolution as in \eqref{eq:res-nullcorr}
\[
0 \lra S_A(-3) \stackrel{P_{2,\epsilon}}{\lra} S_A(-2)^{\oplus 4} \stackrel{P_{1,\epsilon}}{\lra} S_A(-1)^{\oplus 5} \lra M_\epsilon \lra 0
\]
Arguing as in the paragraph after \eqref{eq:res-nullcorr} we conclude that, up to a change of coordinates, $P_{2,\epsilon} = P_2$ and $P_{1,\epsilon}P_2 = 0$ implies that $P_{1,\epsilon} = P_1 + \epsilon Q$ where $QP_2 = 0$. The rows of $Q$ are linear Koszul syzygies of $P_2$. Comparing entries with \eqref{eq:matB} we get a change of coordinates and basis change over $A$ that preserves $P_2$ and sends $P_{1,\epsilon}$ to $P_1$. Thus, $M_\epsilon \cong M_0 \otimes A$. What remains dependent on $\epsilon$ is the morphism $\rho \colon M_\epsilon \to T_S\otimes A$.

After the above change of coordinates we get $\widehat{\Pi}_\epsilon = (\sigma_\epsilon)_*\Pi_\epsilon$, where $\sigma_\epsilon \equiv \id \mod \epsilon$. As in the proof of Lemma \ref{lem:nullcorr-n>=7}, the map $\rho$ lifts to a matrix $C_\epsilon^TB$ where $C_\epsilon = (c_{ij})$ has constant entries in $A$. Define $D_i = \sum_j c_{ij}\del{x_{j}}$. Then $D_j \equiv \del{x_{4+j}} \mod \epsilon$. Up to multiplying $\Pi_\epsilon$ by an invertible element of $A$ we may assume that
\[
\widehat{\Pi}_\epsilon = \sum_{i<j}b_{ij}D_i\wedge D_j.
\]
As in Lemma \ref{lem:nullcorr-n>=7} we can also show that $D_j(x_i) = 0$ for $i\leq 3$. Therefore, $D_j$ is in the span of $\del{x_4}, \dots, \del{x_n}$. Take a change of coordinates $\xi_\epsilon$ that fixes the first four variables and satisfies $(\xi_\epsilon)_*\del{x_{i+4}} = D_i$ and $\xi_\epsilon \equiv \id \mod \epsilon$. Then $(\xi_\epsilon)_*\Pi = \widehat{\Pi}_\epsilon$. Therefore $(\sigma_\epsilon^{-1}\circ \xi_\epsilon)_*\Pi = \Pi_\epsilon$.
\end{proof}

\subsection{The case $n=7$}\label{subsec:nullcorr-p7}
For $n=7$, any null-correlation foliation is special. Following the proofs of Lemma \ref{lem:nullcorr-n>=7} and Lemma \ref{lem:special-to-general} we get that such a foliation is projectively equivalent to the one defined by 
\[
\Pi_0 = \sum_{1\leq i<j\leq 4}b_{ij}\del{x_{3+i}}\wedge\del{x_{3+j}}.
\]
A direct computation shows that the orbit of $[\Pi_0] \in \PP H^0(\extp^2\tp7)$ has dimension $40$. However, $\dim T_{[\Pi_0]}\Fol(2,2,7) = 44$. Then our previous argument does not work. 

\begin{proposition}\label{prop:nullcor-n=7}
Let $\NC_7$ be the closure of the orbit of $[\Pi_0]$ in $\Fol(2,2,7)$. Then $\Fol(2,2,7)$ has a nonreduced component supported on $\NC_7$.
\end{proposition}

\begin{proof}
We will show that the foliation $\sF$ defined by $\Pi_0$ is not in the closure of any other component. First, note that the kernel of $v \mapsto v\wedge \Pi_0$ on linear vector fields is precisely $H^0(\tF)$. Hence $h^0(\tF) = 0$ defines an open neighborhood of $[\Pi_0]$. Hence we may discard components where the general member satisfies $h^0(\tF) \geq 1$. We are left with the rational foliations. On the other hand, a general rational foliation of degree two has non-Kupka singularities of codimension at most $3$. The non-Kupka locus of $\sF$ is defined by the coefficients of $\Pi_0$, which has codimension 4. By semicontinuity of the dimension, $\sF$ cannot be in the boundary of a rational component. Therefore, the null-correlation foliations form an irreducible component of $\Fol(2,2,7)$. By our previous comment, this component has dimension $40$ but the Zariski tangent space $T_{[\Pi_0]}\Fol(2,2,7)$ has dimension $44$. Then it must be nonreduced.
\end{proof}

\section{Classification results}\label{sec:class-dim2}
In this section we prove the main classification results. By combining the results already proved, we demonstrate Theorem \ref{thm:mainA}, which classifies the two-dimensional case. Then we do a reduction step to deal with higher dimensions. 

\begin{proof}[Proof of Theorem~\ref{thm:mainA}]
Fix a two-dimensional foliation $\sF$ of degree two on $\pn$. If $h^0(\tF(-1)) > 0$, then $\sF$ is a linear pullback of a degree-two foliation by curves on $\p{n-1}$ and $\sF \in \lpb{2}{2}{n}$. If $\tF \cong \opn^{\oplus 2}$, then $\sF$ is given by the action of a two-dimensional Lie algebra; hence $\sF$ is in a logarithmic component or an affine-Lie component. If $h^0(\tF(-1)) =0 $, $h^0(\tF) > 0$ and $\tF \not\cong \opn^{\oplus 2}$, Proposition \ref{prop:sub-deg2-alg}, Proposition \ref{prop:class-om-nilp}, and Proposition \ref{prop:class-om-diag} yield that $\sF$ is in either a logarithmic component or an affine-Lie component. If $h^0(\tF) = 0$
then $\tF$ is stable. Lemma \ref{lem:reduc-stable-case} and Lemma \ref{lem:split-to-rational} show that $\sF$ must be either a rational foliation or a null-correlation foliation ($n\geq 7$ in this case). Finally, Proposition \ref{prop:nullcorr-n>=8} and Proposition \ref{prop:nullcor-n=7} prove that null-correlation foliations form an irreducible component. 
\end{proof}

\subsection{Higher dimensions}
The strategy to prove Theorem \ref{thm:mainC-intro} is to reduce to linear pullbacks of known models. We first show that these are either logarithmic or 

\begin{proposition}\label{prop:ra=r-1=log}
Let $\sF$ be a degree-two foliation of dimension $k>2$ on $\pn$. If $H^0(\tF(-1))=0$, then either 
\begin{itemize}
    \item $\sF$ is in a logarithmic component; or
    \item $k=3$ and $\tF = \pi^*\tp3(-1)$ for some linear projection $\pi \colon \pn \dashrightarrow \p3$.
\end{itemize}
\end{proposition}

\begin{proof}
Let $\sF$ be a foliation of degree $2$ and dimension $k>2$ such that $H^0(\tF(-1)) = 0$. If $H^0(\tF^\vee) = 0$, then Lemma \ref{lem:reduc-stable-case} and Lemma \ref{lem:split-to-rational} imply that either $\sF$ is a rational foliation (in particular, logarithmic), or $k=3$ and $\tF = \pi^*\tp3(-1)$ for some linear projection $\pi \colon \pn \dashrightarrow \p3$. 

Suppose $H^0(\tF^\vee) \neq 0$. Then a global section defines a subdistribution $\sG$ of dimension $k-1$ and degree at most one. Since $H^0(\tF(-1)) = 0$, $\deg(\sG) = 1$ and $\tG$ must be stable. Since $\mu(\tG) >0$ and $\mu(\tF/\tG) =0$, by \cite[Lemma 4.10]{CPaun} $\sG$ is a foliation.

By Theorem~\ref{thm:fols-deg1}, $\sG$ is defined by a rational map $\rho:\pn\dashrightarrow\PP(1^{n-k+1},2)$. Thus $\sF$ is the pullback of a foliation induced by a global vector field. By Remark~\ref{rem:log(1,2)}, $\sF\in\logg{1^{n-k+1},2}{k}{n}$.
\end{proof}

Next we prove Theorem \ref{thm:mainC-intro}. 

\begin{proof}[Proof of Theorem \ref{thm:mainC-intro}]
If $H^0(\tF(-1)) = 0$,  Proposition \ref{prop:ra=r-1=log} gives (1) or (5) with the linear pullback being the identity. If $h^0(\tF(-1)) = h > 0$, then $\sF = \pi^*\sG$ where $\sG$ is a degree-two foliation on $\p{n-h}$ such that $H^0(\tG(-1)) =0$, and $\pi \colon \pn \dashrightarrow \p{n-h}$ is a linear projection.

If $\dim(\sG) = 1$ then $\sF \in \lpb{2}{k}{n}$, hence (3). If $\dim(\sG) = 2$ then the proof of Theorem \ref{thm:mainA} gives (1), (2), or (4). If $\dim(\sG) > 2$, Proposition \ref{prop:ra=r-1=log} gives (1) or (5).
\end{proof}

\bibliographystyle{abbrvurl}
\bibliography{bibliography}

@book{OkonekSchneiderSpindler1980,
  author    = {Okonek, Christian and Schneider, Michael and Spindler, Heinz},
  title     = {Vector Bundles on Complex Projective Spaces},
  year      = {1980},
  publisher = {Birkh{\"a}user},
  address   = {Basel}
}

@article {Omegar,
    AUTHOR = {Calvo-Andrade, Omegar},
     TITLE = {Irreducible components of the space of holomorphic foliations},
   JOURNAL = {Math. Ann.},
  FJOURNAL = {Mathematische Annalen},
    VOLUME = {299},
      YEAR = {1994},
    NUMBER = {4},
     PAGES = {751--767},
      ISSN = {0025-5831,1432-1807},
   MRCLASS = {32L30},
  MRNUMBER = {1286897},
MRREVIEWER = {Dominique\ Cerveau},
       DOI = {10.1007/BF01459811},
       URL = {https://doi.org/10.1007/BF01459811},
}

@article {MOJ,
    AUTHOR = {Calvo-Andrade, Omegar and Corr\^{e}a, Maur\'{\i}cio and
              Jardim, Marcos},
     TITLE = {Codimension one holomorphic distributions on the projective
              three-space},
   JOURNAL = {Int. Math. Res. Not. IMRN},
  FJOURNAL = {International Mathematics Research Notices. IMRN},
      YEAR = {2020},
    NUMBER = {23},
     PAGES = {9011--9074},
      ISSN = {1073-7928,1687-0247},
   MRCLASS = {32S65 (14F06)},
  MRNUMBER = {4184615},
MRREVIEWER = {Calum\ Spicer},
       DOI = {10.1093/imrn/rny251},
       URL = {https://doi.org/10.1093/imrn/rny251},
}

@article {CPaun,
    AUTHOR = {Campana, Fr\'{e}d\'{e}ric and P\u{a}un, Mihai},
     TITLE = {Foliations with positive slopes and birational stability of
              orbifold cotangent bundles},
   JOURNAL = {Publ. Math. Inst. Hautes \'{E}tudes Sci.},
  FJOURNAL = {Publications Math\'{e}matiques. Institut de Hautes \'{E}tudes
              Scientifiques},
    VOLUME = {129},
      YEAR = {2019},
     PAGES = {1--49},
      ISSN = {0073-8301,1618-1913},
   MRCLASS = {37F75 (14E30 14M22 32S65)},
  MRNUMBER = {3949026},
MRREVIEWER = {Carla\ Novelli},
       DOI = {10.1007/s10240-019-00105-w},
       URL = {https://doi.org/10.1007/s10240-019-00105-w},
}

@article {CP1,
    AUTHOR = {Cukierman, Fernando and Pereira, Jorge Vit\'{o}rio},
     TITLE = {Stability of holomorphic foliations with split tangent sheaf},
   JOURNAL = {Amer. J. Math.},
  FJOURNAL = {American Journal of Mathematics},
    VOLUME = {130},
      YEAR = {2008},
    NUMBER = {2},
     PAGES = {413--439},
      ISSN = {0002-9327,1080-6377},
   MRCLASS = {32S65},
  MRNUMBER = {2405162},
MRREVIEWER = {M.\ G.\ Soares},
       DOI = {10.1353/ajm.2008.0011},
       URL = {https://doi.org/10.1353/ajm.2008.0011},
}

@article {CPV1,
    AUTHOR = {Cukierman, F. and Pereira, J. V. and Vainsencher, I.},
     TITLE = {Stability of foliations induced by rational maps},
   JOURNAL = {Ann. Fac. Sci. Toulouse Math. (6)},
  FJOURNAL = {Annales de la Facult\'{e} des Sciences de Toulouse.
              Math\'{e}matiques. S\'{e}rie 6},
    VOLUME = {18},
      YEAR = {2009},
    NUMBER = {4},
     PAGES = {685--715},
      ISSN = {0240-2963,2258-7519},
   MRCLASS = {32S65 (37F75)},
  MRNUMBER = {2590385},
MRREVIEWER = {Yohann\ Genzmer},
       URL = {http://afst.cedram.org/item?id=AFST_2009_6_18_4_685_0},
}

@article {CLN,
    AUTHOR = {Cerveau, D. and Lins Neto, A.},
     TITLE = {Irreducible components of the space of holomorphic foliations
              of degree two in {$\mathbb C{\rm P}(n)$}, {$n\geq 3$}},
   JOURNAL = {Ann. of Math. (2)},
  FJOURNAL = {Annals of Mathematics. Second Series},
    VOLUME = {143},
      YEAR = {1996},
    NUMBER = {3},
     PAGES = {577--612},
      ISSN = {0003-486X,1939-8980},
   MRCLASS = {32L30 (32G99)},
  MRNUMBER = {1394970},
MRREVIEWER = {Alain\ H\'{e}naut},
       DOI = {10.2307/2118537},
       URL = {https://doi.org/10.2307/2118537},
}

@article {CLNdeg3,
    AUTHOR = {Cerveau, Dominique and Lins Neto, Alcide},
     TITLE = {A structural theorem for codimension-one foliations on {$\mathbb
              P^n$}, {$n\geq 3$}, with an application to degree-three
              foliations},
   JOURNAL = {Ann. Sc. Norm. Super. Pisa Cl. Sci. (5)},
  FJOURNAL = {Annali della Scuola Normale Superiore di Pisa. Classe di
              Scienze. Serie V},
    VOLUME = {12},
      YEAR = {2013},
    NUMBER = {1},
     PAGES = {1--41},
      ISSN = {0391-173X,2036-2145},
   MRCLASS = {37F75 (34M45)},
  MRNUMBER = {3088436},
MRREVIEWER = {Arseniy\ A.\ Shcherbakov},
}

@article {CLNlog,
    AUTHOR = {Cerveau, Dominique and Lins Neto, Alcides},
     TITLE = {Logarithmic foliations},
   JOURNAL = {Ann. Fac. Sci. Toulouse Math. (6)},
  FJOURNAL = {Annales de la Facult\'{e} des Sciences de Toulouse.
              Math\'{e}matiques. S\'{e}rie 6},
    VOLUME = {30},
      YEAR = {2021},
    NUMBER = {3},
     PAGES = {561--618},
      ISSN = {0240-2963,2258-7519},
   MRCLASS = {37F75 (32S65 34M15)},
  MRNUMBER = {4331310},
MRREVIEWER = {Arseniy\ A.\ Shcherbakov},
       DOI = {10.5802/afst.1685},
       URL = {https://doi.org/10.5802/afst.1685},
}

@book{CD2,
 author = {D{\'e}serti, Julie and Cerveau, Dominique},
 title = {Feuilletages et actions de groupes sur les espaces projectifs},
 fseries = {M{\'e}moires de la Soci{\'e}t{\'e} Math{\'e}matique de France. Nouvelle S{\'e}rie},
 series = {M{\'e}m. Soc. Math. Fr., Nouv. S{\'e}r.},
 issn = {0249-633X},
 volume = {103},
 year = {2005},
 publisher = {Soci{\'e}t{\'e} Math{\'e}matique de France (SMF), Paris},
 language = {French},
 doi = {10.24033/msmf.415},
 zbMATH = {5023912},
 Zbl = {1107.37037}
}

@misc{AVdeg4,
      title={Holomorphic foliations of degree four on the complex projective space}, 
      author={Arturo Fernández-Pérez and Vângellis Sagnori Maia},
      year={2022},
      eprint={2208.04092},
      archivePrefix={arXiv},
      primaryClass={math.CV}
}

@article {CLPdeg3,
    AUTHOR = {da Costa, Raphael Constant and Lizarbe, Ruben and Pereira,
              Jorge Vit\'{o}rio},
     TITLE = {Codimension one foliations of degree three on projective
              spaces},
   JOURNAL = {Bull. Sci. Math.},
  FJOURNAL = {Bulletin des Sciences Math\'{e}matiques},
    VOLUME = {174},
      YEAR = {2022},
     PAGES = {Paper No. 103092, 39},
      ISSN = {0007-4497},
   MRCLASS = {37F75 (32G13)},
  MRNUMBER = {4354288},
MRREVIEWER = {F.\ Touzet},
       DOI = {10.1016/j.bulsci.2021.103092},
       URL = {https://doi.org/10.1016/j.bulsci.2021.103092},
}

@article {Medeiros,
    AUTHOR = {de Medeiros, Airton S.},
     TITLE = {Singular foliations and differential {$p$}-forms},
   JOURNAL = {Ann. Fac. Sci. Toulouse Math. (6)},
  FJOURNAL = {Annales de la Facult\'{e} des Sciences de Toulouse.
              Math\'{e}matiques. S\'{e}rie 6},
    VOLUME = {9},
      YEAR = {2000},
    NUMBER = {3},
     PAGES = {451--466},
      ISSN = {0240-2963,2258-7519},
   MRCLASS = {32S65},
  MRNUMBER = {1842027},
MRREVIEWER = {M.\ G.\ Soares},
       URL = {http://www.numdam.org/item?id=AFST_2000_6_9_3_451_0},
}

@article {GJM,
    AUTHOR = {Galeano, Hugo and Jardim, Marcos and Muniz, Alan},
     TITLE = {Codimension one distributions of degree 2 on the
              three-dimensional projective space},
   JOURNAL = {J. Pure Appl. Algebra},
  FJOURNAL = {Journal of Pure and Applied Algebra},
    VOLUME = {226},
      YEAR = {2022},
    NUMBER = {2},
     PAGES = {Paper No. 106840, 32},
      ISSN = {0022-4049,1873-1376},
   MRCLASS = {58A17 (13D02 14D20 14D22 14J60)},
  MRNUMBER = {4287792},
MRREVIEWER = {Hugo\ Torres\ L\'{o}pez},
       DOI = {10.1016/j.jpaa.2021.106840},
       URL = {https://doi.org/10.1016/j.jpaa.2021.106840},
}

@book {Jou,
    AUTHOR = {Jouanolou, J. P.},
     TITLE = {\'{E}quations de {P}faff alg\'{e}briques},
    SERIES = {},
    VOLUME = {708.},
 PUBLISHER = {Springer, Berlin},
      YEAR = {1979},
     PAGES = {v+255},
      ISBN = {3-540-09239-0},
   MRCLASS = {14D05 (32C40 34A20 35A20 58A17)},
  MRNUMBER = {537038},
MRREVIEWER = {V.\ A.\ Golubeva},
}

@article {LPT2,
    AUTHOR = {Loray, Frank and Pereira, Jorge Vit\'{o}rio and Touzet,
              Fr\'{e}d\'{e}ric},
     TITLE = {Foliations with trivial canonical bundle on {F}ano 3-folds},
   JOURNAL = {Math. Nachr.},
  FJOURNAL = {Mathematische Nachrichten},
    VOLUME = {286},
      YEAR = {2013},
    NUMBER = {8-9},
     PAGES = {921--940},
      ISSN = {0025-584X,1522-2616},
   MRCLASS = {32S65 (14J45 37F75 53D17)},
  MRNUMBER = {3066408},
MRREVIEWER = {Ahmed\ Lesfari},
       DOI = {10.1002/mana.201100354},
       URL = {https://doi.org/10.1002/mana.201100354},
}

@article {LPT3,
    AUTHOR = {Loray, Frank and Pereira, Jorge Vitorio and Touzet,
              Fr\'{e}d\'{e}ric},
     TITLE = {Deformation of rational curves along foliations},
   JOURNAL = {Ann. Sc. Norm. Super. Pisa Cl. Sci. (5)},
  FJOURNAL = {Annali della Scuola Normale Superiore di Pisa. Classe di
              Scienze. Serie V},
    VOLUME = {21},
      YEAR = {2020},
     PAGES = {1315--1331},
      ISSN = {0391-173X,2036-2145},
   MRCLASS = {37F75 (32G07 32S65)},
  MRNUMBER = {4288634},
MRREVIEWER = {Calum\ Spicer},
       DOI = {10.2422/2036-2145.201712\_010},
       URL = {https://doi.org/10.2422/2036-2145.201712_010},
}

@misc{C,
      title={Splitting aspects of holomorphic distributions with locally free tangent sheaf}, 
      author={Raphael Constant da Costa},
      year={2024},
      eprint={2405.17415},
      archivePrefix={arXiv},
      primaryClass={math.CV}
}

@Article{H2,
 Author = {Robin {Hartshorne}},
 Title = {{Stable reflexive sheaves}},
 FJournal = {{Mathematische Annalen}},
 Journal = {{Math. Ann.}},
 ISSN = {0025-5831; 1432-1807/e},
 Volume = {254},
 Pages = {121--176},
 Year = {1980},
 Publisher = {Springer, Berlin/Heidelberg},
 Language = {English},
 MSC2010 = {14F05 32L20 14D22 57R20},
 Zbl = {0431.14004},
 DOI= {10.1007/BF01467074}
}

@misc{PdS,
      title={Distributions with locally free tangent sheaf}, 
      author={J. V. Pereira and J. P. dos Santos},
      year={2024},
      eprint={2408.10057},
      archivePrefix={arXiv},
      primaryClass={math.AG},
      url={https://arxiv.org/abs/2408.10057}, 
}

@book{Borel,
 author = {Borel, Armand},
 title = {Linear algebraic groups.},
 edition = {2nd enlarged ed.},
 fseries = {Graduate Texts in Mathematics},
 series = {Grad. Texts Math.},
 issn = {0072-5285},
 volume = {126},
 isbn = {0-387-97370-2},
 year = {1991},
 publisher = {New York etc.: Springer-Verlag},
 language = {English},
 zbMATH = {50185},
 Zbl = {0726.20030}
}

@book{HL-sheaves,
 author = {Huybrechts, Daniel and Lehn, Manfred},
 title = {The geometry of moduli spaces of sheaves},
 edition = {2nd ed.},
 isbn = {978-0-521-13420-0},
 year = {2010},
 publisher = {Cambridge: Cambridge University Press},
 language = {English},
 zbMATH = {5765056},
 Zbl = {1206.14027}
}

@book{Matsumura,
 author = {Matsumura, Hideyuki},
 title = {Commutative ring theory. {Transl}. from the {Japanese} by {M}. {Reid}.},
 edition = {Paperback ed.},
 fseries = {Cambridge Studies in Advanced Mathematics},
 series = {Camb. Stud. Adv. Math.},
 volume = {8},
 isbn = {0-521-36764-6},
 year = {1989},
 publisher = {Cambridge etc.: Cambridge University Press},
 language = {English},
 zbMATH = {43569},
 Zbl = {0666.13002}
}

@book{BH-CMrings,
 author = {Bruns, Winfried and Herzog, J{\"u}rgen},
 title = {Cohen-Macaulay rings.},
 edition = {Rev. ed.},
 fseries = {Cambridge Studies in Advanced Mathematics},
 series = {Camb. Stud. Adv. Math.},
 volume = {39},
 isbn = {0-521-56674-6},
 year = {1998},
 publisher = {Cambridge: Cambridge University Press},
 language = {English},
 zbMATH = {1194481},
 Zbl = {0909.13005}
}

@article{HKS,
 author = {Herzog, J{\"u}rgen and Kumashiro, Shinya and Stamate, Dumitru I.},
 title = {Graded {Bourbaki} ideals of graded modules},
 fjournal = {Mathematische Zeitschrift},
 journal = {Math. Z.},
 issn = {0025-5874},
 volume = {299},
 number = {3-4},
 pages = {1303--1330},
 year = {2021},
 language = {English},
 doi = {10.1007/s00209-021-02724-8},
 zbMATH = {7428860},
 Zbl = {1483.13006}
}

@book{Gantm.vol1,
 author = {Gantmacher, F. R.},
 title = {The theory of matrices. {Vol}. 1. {Transl}. from the {Russian} by {K}. {A}. {Hirsch}.},
 edition = {Reprint of the 1959 translation},
 isbn = {0-8218-1376-5; 0-8218-1393-5},
 year = {1998},
 publisher = {Providence, RI: AMS Chelsea Publishing},
 language = {English},
 zbMATH = {1235881},
 Zbl = {0927.15001}
}

@article{BeRo,
 author = {Bertone, Cristina and Roggero, Margherita},
 title = {Positivity of {Chern} classes for reflexive sheaves on {{\(\mathbb{P}^{N }\)}}},
 fjournal = {Geometriae Dedicata},
 journal = {Geom. Dedicata},
 issn = {0046-5755},
 volume = {142},
 pages = {121--138},
 year = {2009},
 language = {English},
 doi = {10.1007/s10711-009-9362-5},
 zbMATH = {5614992},
 Zbl = {1174.14015}
}

@book{DN-poisson,
 author = {Dufour, Jean-Paul and Nguyen Tien Zung},
 title = {Poisson structures and their normal forms},
 fseries = {Progress in Mathematics},
 series = {Prog. Math.},
 issn = {0743-1643},
 volume = {242},
 isbn = {3-7643-7334-2},
 year = {2005},
 publisher = {Basel: Birkh{\"a}user},
 language = {English},
 zbMATH = {2237399},
 Zbl = {1082.53078}
}

@misc{stacks-project,
  author       = {The {Stacks project authors}},
  title        = {The Stacks project},
  howpublished = {\url{https://stacks.math.columbia.edu}},
  year         = {2026},
}

@book{Eisen,
 author = {Eisenbud, David},
 title = {Commutative algebra. {With} a view toward algebraic geometry},
 fseries = {Graduate Texts in Mathematics},
 series = {Grad. Texts Math.},
 issn = {0072-5285},
 volume = {150},
 isbn = {3-540-94269-6; 3-540-94268-8},
 year = {1995},
 publisher = {Berlin: Springer-Verlag},
 language = {English},
 zbMATH = {704831},
 Zbl = {0819.13001}
}

@article{CLNpullback,
 author = {Camacho, Cesar and Lins Neto, Alcid{\'e}s},
 title = {The topology on integrable differential forms near a singularity},
 fjournal = {Publications Math{\'e}matiques},
 journal = {Publ. Math., Inst. Hautes {\'E}tud. Sci.},
 issn = {0073-8301},
 volume = {55},
 pages = {5--35},
 year = {1982},
 language = {English},
 doi = {10.1007/BF02698693},
 url = {https://eudml.org/doc/103982},
 zbMATH = {3796868},
 Zbl = {0505.58026}
}
\end{document}